\documentclass[fleqn,a4paper,10pt]{amsart}

\usepackage{graphicx}
\usepackage{tikz}
\usetikzlibrary{shapes,decorations.pathmorphing}
\usepackage{amsmath}
\usepackage{amssymb,amsthm}
\usepackage{comment}
\usepackage{hyperref}

\theoremstyle{definition}

\title[Refined hull laws]{Refined universal laws for hull volumes and perimeters in large planar maps}
\author{Emmanuel Guitter}
\address{Institut de physique th\'eorique, Universit\'e Paris Saclay, CEA, CNRS, F-91191 Gif-sur-Yvette}
\email{emmanuel.guitter@cea.fr}

\begin{document}
\maketitle

\begin{abstract}
We consider ensembles of planar maps with two marked vertices at distance $k$ from each other and 
look at the closed line separating these vertices and lying at distance $d$ from the first one ($d<k$).
This line divides the map into two components, the hull at distance $d$ which corresponds to the
part of the map lying on the same side as the first vertex and its complementary. The number of faces 
within the hull is called the hull volume and the length of the separating line the hull perimeter. We study
the statistics of the hull volume and perimeter for arbitrary $d$ and $k$ in the limit of infinitely large planar
quadrangulations, triangulations and Eulerian triangulations. We consider more precisely situations where
both $d$ and $k$ become large with the ratio $d/k$ remaining finite. For infinitely large maps, 
two regimes may be encountered: either the hull has a finite volume and its complementary is infinitely large, or the
hull itself has an infinite volume and its complementary is of finite size.
We compute the probability for the map to be in either regime as a function of $d/k$ as well as a number of 
universal statistical laws for the hull perimeter and volume when maps are conditioned to be in one regime or the other.
\end{abstract}

\section{Introduction}
\label{sec:introduction}

The study of random planar maps, which are connected graphs embedded on the sphere, has been for more than fifty years the subject of some intense
activity among combinatorialists and probabilists, as well as among physicists in various domains. Very recently, some special attention was
paid to statistical properties of the \emph{hull} in random planar maps, a problem which may be stated as follows: consider an ensemble of planar maps having
two marked vertices at graph distance $k$ from each other. For any non-negative $d$ strictly less than $k$, we may find  
a closed line ``at distance $d$" (i.e.\ made of edges connecting vertices at distance $d$ or so) from the first vertex and separating the two marked vertices 
from each other. Several prescriptions may be adopted for a univocal definition of this separating line but they all eventually give rise to similar statistical properties. 
The separating line divides de facto the map into two connected components, each containing one of the marked vertices. The \emph{hull at distance $d$} corresponds to the
part of the map lying on the same side as the first vertex (i.e.\ that from which distances are measured). The geometrical characteristics of this hull for arbitrary $d$ and $k$ provide
random variables whose statistics may be studied by various techniques. In particular, the statistics of the \emph{volume of the hull}.
which is its number of faces, and of the \emph{hull perimeter}, which is the length of the separating line, have been the subject of several 
investigations \cite{Krikun03,Krikun05,CLG14a,CLG14b,G16a,Men16,G16b}.

In a recent paper \cite{G16a}, we presented a number of results on the statistics of the hull perimeter at distance $d$ for planar triangulations 
(maps with faces of degree three) and quadrangulations (faces of degree four) in a universal regime of infinitely large maps where both $d$ and $k$ are 
large and \emph{remain of the same order}
(i.e.\ the ratio $d/k$ is kept fixed). As we shall see, for such a regime, although the hull perimeter remains finite (but large, of order $d^2$), the volume of the hull 
at distance $d$ may very well be itself strictly infinite. We will compute below the probability for this to happen, a probability which
remains non-zero for large $d$ and $k$ (unless $d/k\to 0$). In particular, if we wish a non-trivial description of the
hull volume statistics, we have to \emph{condition the map configurations so that their hull volume remains finite}. More generally, 
we may reconsider the statistics of the hull perimeter by separating the contribution coming from the set of map configurations with a finite hull volume
from that coming from the set of map configurations with an infinite hull volume. More simply, we may consider the \emph{hull perimeter conditional statistics} obtained
by limiting the configurations to either set of configurations. It is the subject of the present paper to give a precise description of this
\emph{refined hull statistics} where we control the finite or infinite nature of the hull volume. 
Most of the obtained laws crucially depend on the value of $d/k$ but are the same for planar triangulations and planar quadrangulations, 
as well as for Eulerian triangulations (maps with alternating black and white triangular faces).   

The paper is organized as follows: we first present in Section~\ref{sec:summary} a summary of our results  and give explicit
expressions for the probability to have a finite or an infinite hull volume as a function of $d/k$ (Sect.~\ref{sec:inandout}), as well as
for the conditional probability density for the hull perimeter in both situations (Sect.~\ref{sec:pd}). We then give (Sect.~\ref{sec:jointlaw}) the joint law for the 
hull perimeter and hull volume, assuming that the latter is finite. Section~\ref{sec:strategy} presents the strategy that we use for our calculations
which is based on already known generating functions whose expressions are recalled in the case of quadrangulations (Sect.~\ref{sec:genfunc}). 
We explain in details (Sect.~\ref{sec:extract} ) how to extract from these generating functions the desired statistical results. This strategy is implemented for quadrangulations
in Section~\ref{sec:explicit} where we compute the probability to have a finite or an infinite hull volume (Sect.\ref{sec:explprobs}), the probability density for the hull perimeter
in both regimes (Sect.~\ref{sec:explperimlaw}) and the joint law for the hull perimeter and volume when the latter is finite (Sect.~\ref{sec:expljointlaw}).
Section~\ref{sec:Other} briefly discusses triangulations and Eulerian triangulations for which the same universal laws
as those found in the previous sections for quadrangulations are recovered. We gather a few concluding remarks in Section~\ref{sec:conclusion} and present additional 
non-universal expressions at finite $d$ and $k$ in appendix A.

\section{Summary of the results}
\label{sec:summary}
The results presented in this paper have been obtained for three families of planar maps: (i) planar quadrangulations, i.e.\ 
planar maps whose all faces have degree four, (ii) planar triangulations, i.e.\ planar maps whose all faces have degree three
and (iii) planar Eulerian triangulations,  which are planar triangulations whose faces are colored in black and white with adjacent
faces being of different color. For all these families, we obtain in the limit of large maps \emph{the same laws} for hull volumes and perimeters, up to two non-universal
normalization factors, one for the volume and one for the perimeter (called $f$ and $c$ respectively). The hull volumes and perimeters are defined
as follows: for the three families of maps and for some integer $k\geq 1$, we consider more precisely $k$-pointed-rooted maps, i.e.\ maps 
with a marked vertex $x_0$ (called the origin) and a marked oriented edge pointing from a vertex $x_1$ at graph distance $k$ from the origin $x_0$ 
to a neighbor of $x_1$ at distance $k\!-\!1$ (such neighbor always exists)\footnote{In case (iii), we use more precisely some natural ``oriented graph distance"
using oriented paths keeping black faces on their left, see \cite{G16b}.}. Given $k \geq 3$ and some integer $d$ in the range $2\leq d\leq k\!-\!1$, there exists a simple
closed line along edges of the map, ``at distance $d$" from the origin\footnote{In practice, the line
may be chosen in case (ii) so as to visit only vertices at distance $d$, and in case (i) and (iii) so as to visit alternately vertices at distance $d$ and $d-1$.}
which separates the origin $x_0$ from $x_1$. Several prescription are possible for a univocal definition of this \emph{separating line} and we will adopt here that proposed in
\cite{G16a} in cases (i) and (ii) and in \cite{G16b} for case (iii). We expect that other choices should not modify our results, except possibly for the value 
of the perimeter normalization factor $c$.
The \emph{hull at distance $d$ in our $k$-pointed-rooted map} is defined as the domain of the map lying \emph{on the same side as 
the origin} of the separating line at distance $d$. Its \emph{volume ${\mathcal V}(d)$} is its number of faces and its \emph{perimeter ${\mathcal L}(d)$} the length (i.e.\ number of edges) of 
its boundary, namely the length of the separating line at distance $d$ itself. 

\begin{figure}
\begin{center}
\includegraphics[width=11cm]{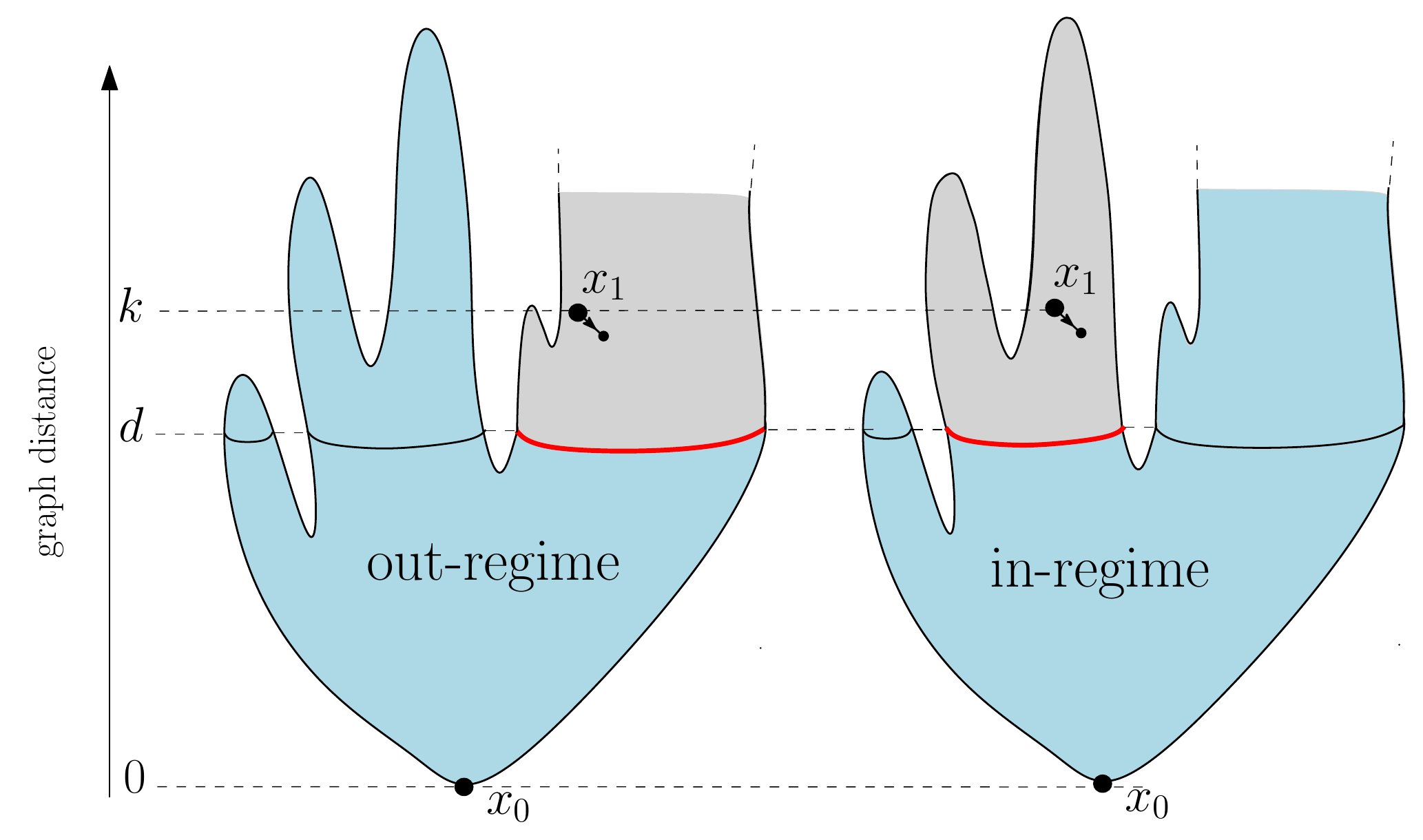}
\end{center}
\caption{An illustration of map configurations in the out- and in-regimes. The map is represented schematically with vertices placed at a height equal to their distance from 
the origin $x_0$. The vertices at distance $d$ from $x_0$ form a number of closed curves at height $d$, one of which (in red) separates $x_0$ from $x_1$ and 
defines the separating line at distance $d$. The part of the map lying on the same side of this separating line as $x_0$ constitutes the hull at distance $d$ (here in light blue). 
For maps with an infinite volume, the configuration is in the out-regime if the hull volume remains finite (configuration on the left) or it is in the in-regime
if the hull volume itself becomes infinite (configuration on the right).}
\label{fig:hullboundary}
\end{figure}
\subsection{The out- and in-regimes}
\label{sec:inandout}
Our results deal with the statistics of uniformly drawn $k$-pointed-rooted maps in the families (i), (ii) or (iii) having a \emph{fixed number of faces} $N$, 
and for a \emph{fixed value of the parameter $k$} and, more precisely, with \emph{the limit $N\to \infty$ of this ensemble}, keeping $k$ finite. 
This corresponds to the so called \emph{local limit} of infinitely large maps and, as in \cite{G16a}, 
we shall denote by $P_k(\{\cdot\})$ the probability of some event $\{\cdot\}$ and $E_k(\{\cdot\})$ the expectation value of some quantity 
$\{\cdot\}$ in this limit. 

In the limit $N\to \infty$, two situations may occur: either the volume $\mathcal V(d)$ remains finite and the number of faces $N\!-\!{\mathcal V}(d)$ of
the complementary of the hull (namely the part of the map lying on the same side of the separating line as $x_1$) is infinite.
This situation will be referred to as the \emph{``out-regime"} in the following. Or the volume ${\mathcal V}(d)$ is itself infinite while the number $N\!-\!{\mathcal V}(d)$ 
of faces of the complementary of the hull remains finite. This situation will be referred to as the \emph{``in-regime"} in the following.
The case where both ${\mathcal V}(d)$ and $N\!-\!{\mathcal V}(d)$ would be infinite is expected to be suppressed when $N\to \infty$
(i.e.\ the number of configurations in this regime does not grow with $N$ as fast as that in the out- and in-regimes).  Situations in the
out- and in- regime are illustrated in figure \ref{fig:hullboundary}. 

The main novelty in this paper is that our laws will discriminate between situations where the map configurations are in the out- or in the in-regime. 
We shall use accordingly the notations
\begin{equation*}
\begin{split}
&P^{\rm{out}}_{k,d}(\{\cdot\})=P_k\big(\{\cdot\}\ \hbox{and}\ {\mathcal V}(d)\ \hbox{finite}\big)\ ,
\quad P^{\rm{in}}_{k,d}(\{\cdot\})=P_k\big(\{\cdot\}\ \hbox{and}\ {\mathcal V}(d)\ \hbox{infinite}\big)\ ,\\
& E^{\rm{out}}_{k,d}(\{\cdot\})=E_k\big(\{\cdot\}\!\times\!\theta_{\rm finite}({\mathcal V}(d)\big)\ ,
\quad E^{\rm{in}}_{k,d}(\{\cdot\})=E_k\big(\{\cdot\}\!\times\!(1-\theta_{\rm finite}({\mathcal V}(d))\big)\,\\
\end{split}
\end{equation*}  
with $\theta_{\rm finite}({\mathcal V})=1$ if ${\mathcal V}$ is finite and $0$ otherwise.
Alternatively, we will consider conditional probabilities and conditioned expectation values, defined respectively as
\begin{equation*}
\begin{split}
&P_k(\{\cdot\}|{\mathcal V}(d)\ \hbox{finite})=\frac{P^{\rm{out}}_{k,d}(\{\cdot\})}{P_k({\mathcal V}(d)\ \hbox{finite})}\ ,
\quad P_k(\{\cdot\}|{\mathcal V}(d)\ \hbox{infinite})=\frac{P^{\rm{in}}_{k,d}(\{\cdot\})}{P_k({\mathcal V}(d)\ \hbox{infinite})}\ .\\
&E_k(\{\cdot\}|{\mathcal V}(d)\ \hbox{finite})=\frac{E^{\rm{out}}_{k,d}(\{\cdot\})}{P_k({\mathcal V}(d)\ \hbox{finite})}\ ,
\quad E_k(\{\cdot\}|{\mathcal V}(d)\ \hbox{infinite})=\frac{E^{\rm{in}}_{k,d}(\{\cdot\})}{P_k({\mathcal V}(d)\ \hbox{infinite})}\ ,\\
\end{split}
\end{equation*}  
where $P_k({\mathcal V}(d)\ \hbox{finite})=E_k\big(\theta_{\rm finite}({\mathcal V}(d)\big)=1-P_k({\mathcal V}(d)\ \hbox{infinite})$.
\vskip .3cm
Universal laws may are obtained when $k$ and $d$ themselves become large simultameously, i.e.\ upon taking the limit $k\to \infty$, $d\to \infty$ with
$d/k$ fixed (necessarily between $0$ and $1$). We set accordingly 
\begin{equation*}
u\equiv \frac{d}{k}\ , \qquad 0\leq u\leq 1\ ,
\end{equation*}
and our universal results will deal with configurations having \emph{a fixed value of $u$}. 

As in \cite{G16a}, we insist on that we first let $N\to \infty$, and only then take the limit of large $k$ and $d$. In particular, this is to be contrasted with the so called scaling
limit where $N$, $k$ and $d$ would tend simultaneously to infinity with $k\sim d\sim N^{1/4}$.
Note also that our universal laws describe a broader regime than that explored in most papers so far on the hull statistics \cite{Krikun03,Krikun05,CLG14a,CLG14b,Men16},
where the hull boundary is defined as a closed line \emph{separating some origin vertex $x_0$ from infinity} in pointed maps of infinite size.
This latter, more restricted, regime may be recovered in our framework by sending first $k\to \infty$ with $d$ kept finite, and only then letting eventually $d\to \infty$.
As we shall discuss, the results obtained for this latter order of limits match precisely those obtained by taking the limit $u\to 0$ of our results and they may thus be 
considered as particular instances of our more general laws for arbitrary $u$. To be precise, we \emph{observe} that, for all the observables $\{\cdot\}_{d}$ depending on $d$ that we consider, we have the
equivalence
\begin{equation}
\begin{split}
&\lim_{u\to 0} \Big( \lim_{k\to \infty} P_{k}(\{\cdot\}_{k\, u})\Big)=\lim_{d\to \infty}\Big( \lim_{k\to \infty} P_{k}(\{\cdot\}_{d})\Big)\ ,\\
&\lim_{u\to 0} \Big( \lim_{k\to \infty} E_{k}(\{\cdot\}_{k\, u})\Big)=\lim_{d\to \infty}\Big( \lim_{k\to \infty} E_{k}(\{\cdot\}_{d})\Big)\ .\
\end{split}
\label{eq:utozero}
\end{equation}
This equivalence is not a surprise since the limit $u\to 0$ describes precisely situations where the distance $d$ does not scale with $k$. We have however no 
rigorous argument to state that the above identity (based on an inversion of limits) should hold in all generality for any observable $\{\cdot\}_{d}$\footnote{If the observable depends on both $d$ and $k$, the equivalence clearly cannot be true in general
as seen by taking for instance the expectation of $d^2/(k+d^2)$ equal to $1$ or $0$ according to the order of the limits.} .
\begin{figure}
\begin{center}
\includegraphics[width=6.5cm]{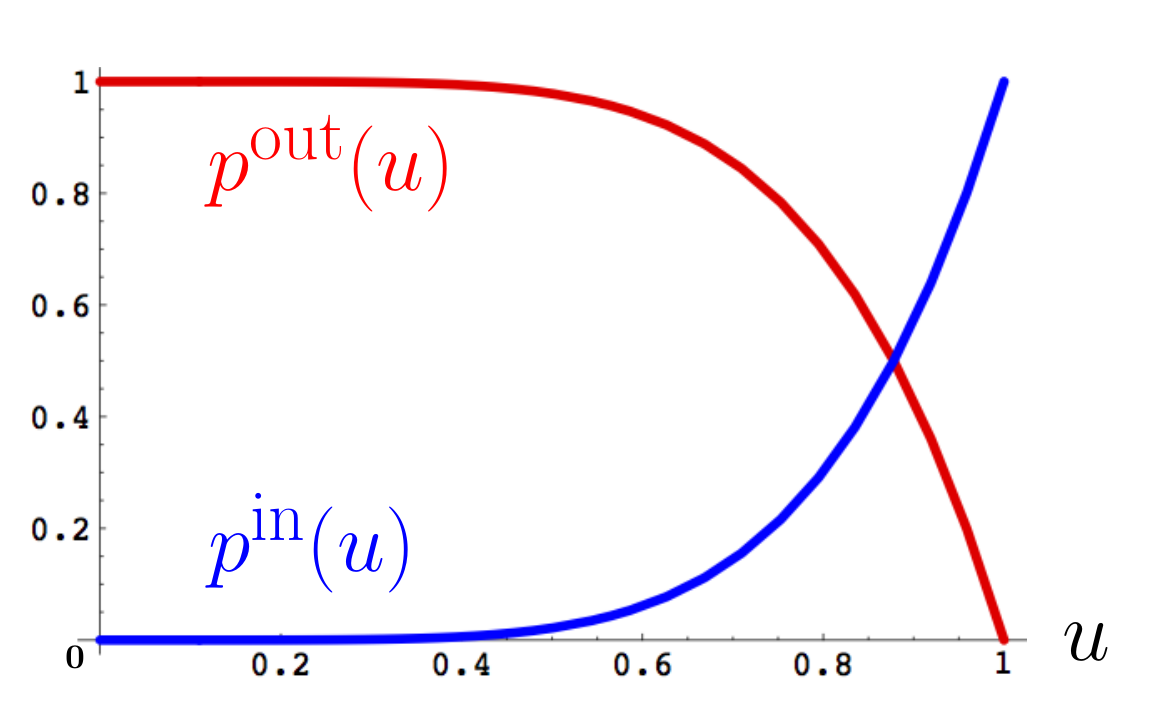}
\end{center}
\caption{A plot of the probability  $p^{\rm out}(u)$ (in red) and the complementary probability $p^{\rm in}(u)$ (in blue)
as a function of $u=d/k$, as given by \eqref{eq:poutin}.}
\label{fig:probasoutin}
\end{figure}
\vskip .3cm
Our first result is an expression, for a given $u$, of the probability that a randomly picked $k$-pointed-rooted map be in the out- or in the in-regime. We find:
\begin{equation}
\begin{split}
&p^{\rm out}(u)\equiv \lim_{k\to \infty} P_k({\mathcal V}(k\, u)\ \hbox{finite})=\frac{1}{4} \left(4-7 u^6+3 u^7\right)\ ,\\
&p^{\rm in}(u)\equiv \lim_{k\to \infty} P_k({\mathcal V}(k\, u)\ \hbox{infinite})=\frac{1}{4} (7-3 u)\, u^6\ ,\\
\end{split}
\label{eq:poutin}
\end{equation}
with of course $p^{\rm out}(u)+p^{\rm in}(u)=1$. Note that these probabilities involve no normalization factor and are the same for the 
three families (i) (ii) and (iii) that we considered. They are represented in figure~\ref{fig:probasoutin}. For $u\to 0$,
we have $p^{\rm out}(u)\to 1$ and $p^{\rm in}(u)\to 0$ so that the map configuration is in the out-regime with probability $1$. This is a first manifestation 
of the equivalence \eqref{eq:utozero} above. Indeed, sending $k\to \infty$ first before letting $d$ be large ensures that the connected component containing
$x_1$ (i.e. the complementary of the hull at distance $d$) has some infinite volume, hence the configuration necessarily lies in the out-regime. 
On the other hand, when $u\to 1$, we see that $p^{\rm out}(u)\to 0$ and $p^{\rm in}(u)\to 1$. This corresponds to situations where 
the vertex $x_1$ remains at a finite distance from the separating line at distance $d$, with $d$ becoming infinitely large.
In such a situation, the connected part containing $x_1$ has a finite volume with probability $1$. This result may be explained 
heuristically as follows: a rough estimate of the probability $p^{\rm out}(1)$ is given by the ratio of the length of the line at distance $d$
separating $x_0$ from infinity by the length of the boundary of the ball of radius $d$ with origin $x_0$. Indeed, the first length measures
the number of ways to place $x_1$ ``just above"\footnote{By ``just above", we mean at a distance from the origin larger than that of the line
by a quantity remaining bounded when $d$ becomes large.} the line separating $x_0$ from infinity while the second length 
is an equivalent measure of the number of ways of placing  $x_1$ anywhere just above a line at distance $d$. Since the first length typically grows like $d^2$ \cite{Krikun03,Krikun05,CLG14a,CLG14b,G16a,Men16,G16b}
while the second length grows like $d^3$ (recall that random maps have fractal dimension $4$ \cite{AmWa95,ChSc04}), the ratio vanishes as $1/d$ when $d\to \infty$
hence $p^{\rm out}(1)$ vanishes.

\subsection{Probability density for the rescaled perimeter in the out- and in-regimes}
\label{sec:pd}
Our second result concerns the probability density for the hull perimeter at distance $d=k\, u$ in the out- and in-regimes. For large $d$,
${\mathcal L}(d)$ scales as $d^2$ so a finite probability density is obtained for the \emph{rescaled perimeter} 
\begin{equation*}
L(d)\equiv \frac{{\mathcal L}(d)}{d^2}\ .
\end{equation*}
We define more precisely the probability densities
\begin{equation*}
\begin{split}
&D^{\rm out}(L,u)\equiv \lim_{k\to \infty} \frac{1}{dL} P^{\rm out}_{k,k\, u}\left(L\leq L(k\, u) < L+dL\right)\ ,\\
&D^{\rm in}(L,u)\equiv  \lim_{k\to \infty} \frac{1}{dL} P^{\rm in}_{k,k\, u}\left(L\leq L(k\, u)  < L+dL\right)\ ,
\end{split}
\end{equation*}
for which we find the following explicit expressions:
\begin{equation}
\begin{split}
&\hskip -0.8cm D^{\rm out}(L,u)=
\frac{(1-u)^4}{2 c \sqrt{\pi } u} \\ &\hskip 0.2cm \times e^{-B X} \left(-2 \sqrt{X} ((X-10) X-2)+e^X \sqrt{\pi } X (X (2 X-5)+6) \left(1-\text{erf}\left(\sqrt{X}\right)\right)\right)\ ,\\
&\hskip -0.8cm D^{\rm in}(L,u)=
\frac{u^5}{2 c \sqrt{\pi } (1-u)^2}\\&\hskip 0.2cm \times 
e^{-B X}  (B X+2) \left(2 \sqrt{X} (X+1)-e^X \sqrt{\pi } X (2 X+3) \left(1-\text{erf}\left(\sqrt{X}\right)\right)\right)\ ,\\
&\hskip -0.8cm \hbox{where}\ X\equiv X(L,u)=\frac{u^2}{(1-u)^2}\frac{L}{c}\ , \qquad B\equiv B(u)=\frac{(1-u)^2}{u^2}\ .\\
\label{eq:probainoutL}
\end{split}
\end{equation}
Here $c$ is a normalization factor given in cases (i), (ii) and (iii) respectively by:
\begin{equation}
\hbox{(i):}\  c=\frac{1}{3}\  , \qquad \hbox{(ii):}\  c=\frac{1}{2}\  , \qquad \hbox{(iii):}\  c=\frac{1}{4}\ .
\label{eq:cval}
\end{equation}
Note that by definition, we have the normalizations
\begin{equation*}
\int_0^\infty D^{\rm out}(L,u)\, dL= p^{\rm out}(u)\ , \qquad \int_0^\infty D^{\rm in}(L,u)\, dL= p^{\rm in}(u)\ ,
\end{equation*}
a result which may be checked directly from the explicit expressions \eqref{eq:poutin} and \eqref{eq:probainoutL}. Note also that the ratio
$D^{\rm out}(L,u)/p^{\rm out}(u)$ (resp. $D^{\rm in}(L,u)/p^{\rm in}(u)$) denotes, at fixed $u=k/d$, the probability
density for the rescaled perimeter ${\mathcal L}(d)/d^2$ for map configurations \emph{conditioned to be in the out-regime} (resp. in the in-regime), with an integral
over $L$ now normalized to $1$.

\begin{figure}
\begin{center}
\includegraphics[width=9cm]{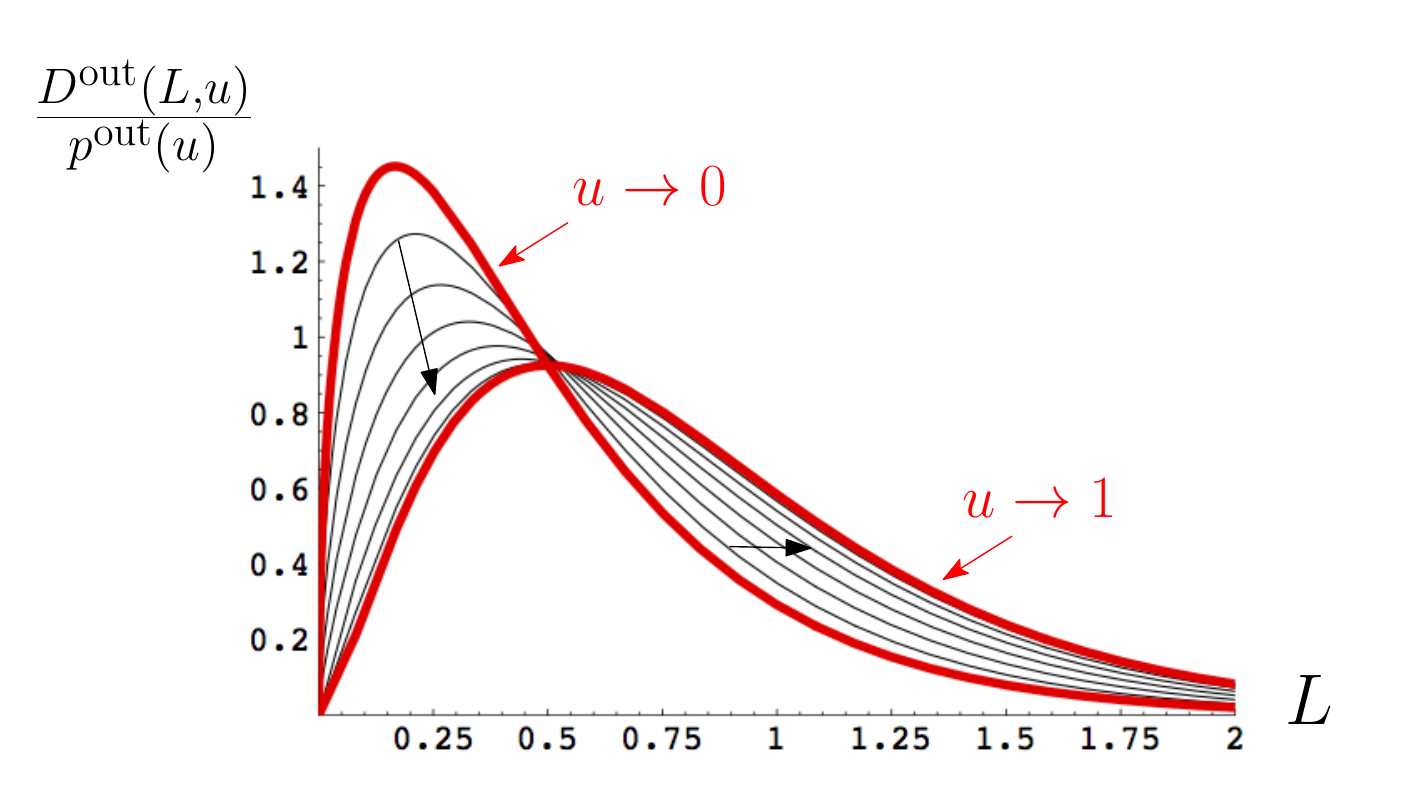}
\end{center}
\caption{The conditional probability density 
$D^{\rm out}(L,u)/p^{\rm out}(u)$ as a function of $L$ for increasing values of $u$ (following the arrow) and its $u\to 0$ and $u\to 1$ limits,
as given by \eqref{eq:uzeroone}.
}
\label{fig:probaoutL}
\end{figure}
\vskip .3cm
Let us now analyze these latter conditional probability densities in more details. Let us first assume that the map configuration lies in the out-regime: the conditional probability density 
$D^{\rm out}(L,u)/p^{\rm out}(u)$, as displayed in figure~\ref{fig:probaoutL}, varies for increasing $u$ between its $u\to 0$ and $u\to 1$ limits
given, from the general expression \eqref{eq:probainoutL}, by
\begin{equation}
\begin{split}
&\lim_{u\to 0} \frac{D^{\rm out}(L,u)}{p^{\rm out}(u)}=2\, \sqrt{L}\, \frac{e^{-\frac{L}{c}}}{c^{3/2} \sqrt{\pi }}\ ,\\
&\lim_{u\to 1} \frac{D^{\rm out}(L,u)}{p^{\rm out}(u)}=\frac{4}{3}(\sqrt{L})^3\frac{e^{-\frac{L}{c}}}{c^{5/2} \sqrt{\pi }}\ .\\
\end{split}
\label{eq:uzeroone}
\end{equation}
It is easy to verify that the $u\to 0$ expression above reproduces precisely the result obtained by first sending $k\to \infty$, and then $d\to\infty$,
in agreement with the announced equivalence \eqref{eq:utozero}. As expected, this expression therefore matches that of Krikun \cite{Krikun03,Krikun05} and of   
Curien and Le Gall \cite{CLG14a,CLG14b} concerning the probability density for the length of the 
line at distance $d$ separating some origin $x_0$ from infinity in large pointed maps of the family at hand (with possibly different values of $c$
due to inequivalent prescriptions for the definition of the separating line). Note that the requirement that the configuration be 
in the out-regime is actually not constraining for $u\to 0$ since $p^{\rm out}(0)=1$.
 
For $u\to 1$, the requirement to be in the out-regime restricts the set of configurations to those where we have
chosen $x_1$ in the vicinity (i.e. just above) the line separating $x_0$ from infinity (so that the domain in which $x_1$ lies is infinite).
As just discussed, this line has a length $L\, d^2$ with density probability $2\, \sqrt{L}\, e^{-\frac{L}{c}}/(c^{3/2} \sqrt{\pi })$ 
while the number of choices for $x_1$ is (for fixed $d$) proportional to $L$. The conditional probability density for $L(d)$ 
in the in-regime is thus expected to be 
\begin{equation*}
\frac{L\times 2\, \sqrt{L}\, \frac{e^{-\frac{L}{c}}}{c^{3/2} \sqrt{\pi }}}{\int_0^\infty L\times 2\, \sqrt{L}\, \frac{e^{-\frac{L}{c}}}{c^{3/2} \sqrt{\pi }}\, dL}
= \frac{4}{3}(\sqrt{L})^3\frac{e^{-\frac{L}{c}}}{c^{5/2} \sqrt{\pi }}\ ,
\end{equation*}
and this is precisely the result obtained above.

\begin{figure}
\begin{center}
\includegraphics[width=10cm]{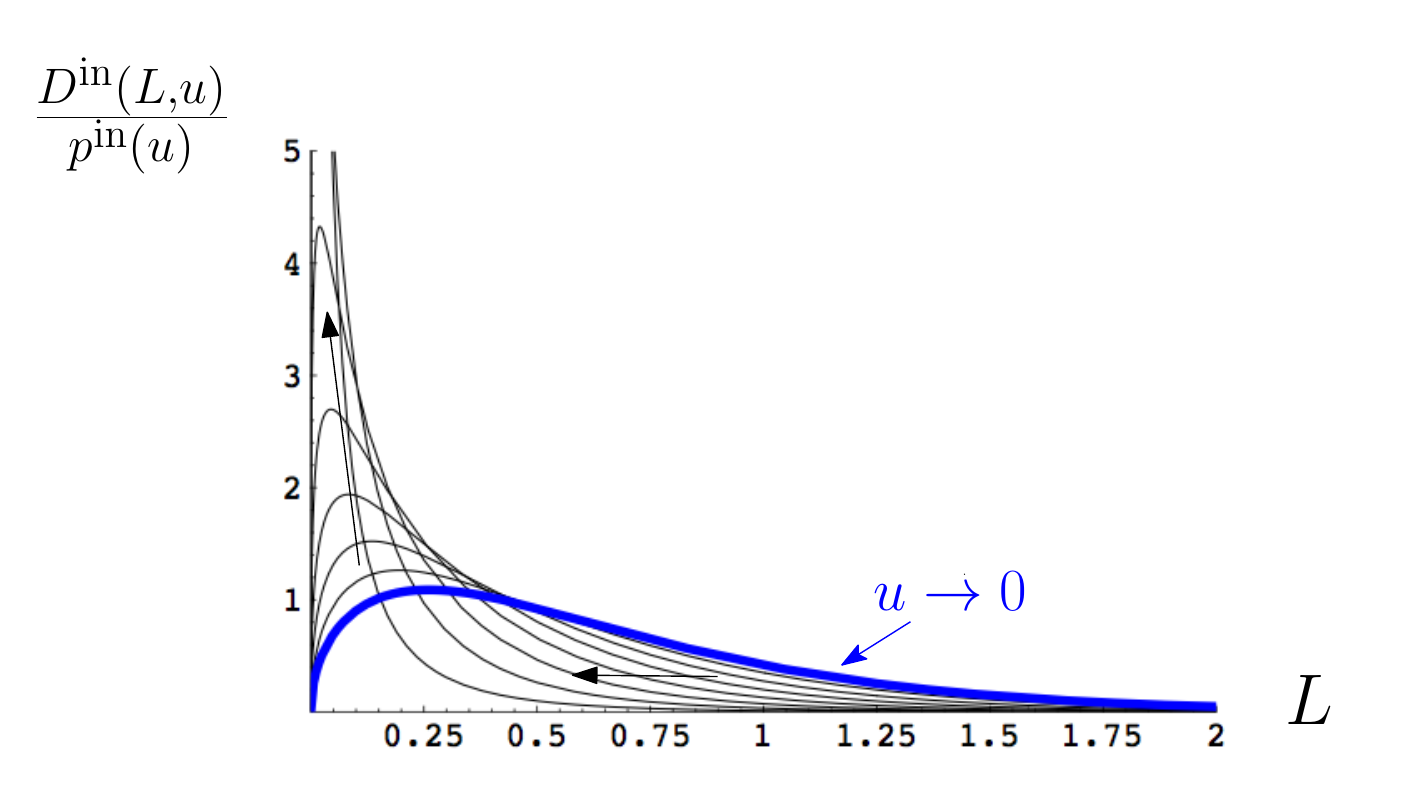}
\end{center}
\caption{The conditional probability density 
$D^{\rm in}(L,u)/p^{\rm in}(u)$ as a function of $L$ for increasing values of $u$ (following the arrow) and its $u\to 0$ limit,
as given by \eqref{eq:uzeroin}}
\label{fig:probainL}
\end{figure}
\begin{figure}
\begin{center}
\includegraphics[width=10cm]{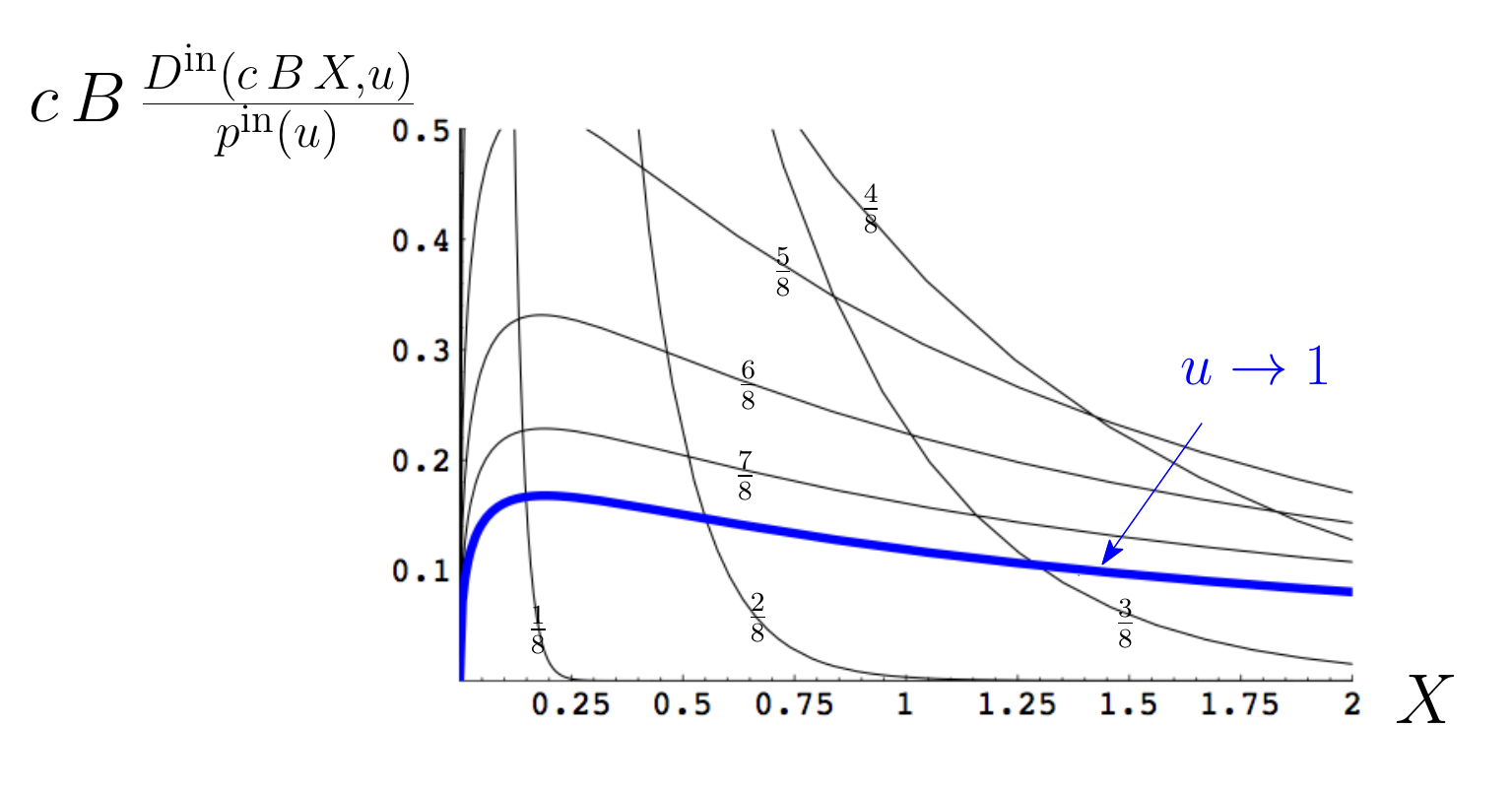}
\end{center}
\caption{The conditional probability density in the in-regime for the variable $X= L/(c B)$, at $u=1/8,2/8,3/8,\cdots$ and in the limit $u\to 1$.}
\label{fig:probainLrescaled}
\end{figure}

Let us now assume that the map configuration lies in the in-regime and discuss the corresponding conditional probability density for $L(d)$. 
As displayed in figure~\ref{fig:probainL}, $D^{\rm in}(L,u)/p^{\rm in}(u)$
varies for increasing $u$ between its $u\to 0$ limit given, from the general expression \eqref{eq:probainoutL}, by
\begin{equation}
\lim_{u\to 0} \frac{D^{\rm in}(L,u)}{p^{\rm in}(u)}=\frac{4}{7}\, \sqrt{L}\,  (2 c+L)\, \frac{e^{-\frac{L}{c}}}{ c^{5/2} \sqrt{\pi }}
\label{eq:uzeroin}
\end{equation}
and a degenerate $u\to 1$ limit where only the rescaled length $L=0$ is selected. Recall that $p^{\rm in}(u)\to 0$ for $u\to 0$ and the limiting law just above
for $u\to 0$ therefore describes a very restricted set of configurations where the connected domain containing $x_1$, although $k$ becomes arbitrary larger
than $d$, remains of finite volume. As for the $u\to 1$ limit, the fact that the probability density concentrates around $L=0$ means that ${\mathcal L}(d)$ 
scales less rapidly than $d^2$ in this limit and that some new appropriate rescaling is required. As already discussed in \cite{G16a},
a non-trivial law is in fact obtained by switching to the variable $X$ in \eqref{eq:probainoutL}, i.e.\ considering the probability density for the rescaled length
\begin{equation*}
X(k,d)\equiv \frac{{\mathcal L}(d)}{c\, (k^2-d^2)}=\frac{u^2}{c\, (1-u)^2}\ L(d) \ ,
\end{equation*}
where the coefficient $c$ is arbitrarily included in the definition of $X(k,d)$ so as to have the same limiting law for the three map families (i), (ii) and (iii).
Setting $B=(1-u)^2/u^2$ as in \eqref{eq:probainoutL}, so that $L(d)=c\ B\, X(k,d)$, the probability density for $X(k,k\, u)$ is given for $u\to 1$ by
\begin{equation*}
\lim_{u\to 1} c\, B\, \frac{D^{\rm in}(c\, B\, X,u)}{p^{\rm in}(u)}=\frac{2 \sqrt{X} (X+1)-e^X \sqrt{\pi } X (2 X+3) \left(1-\text{erf}\left(\sqrt{X}\right)\right)}{\sqrt{\pi }}\ .
\end{equation*}
This result matches that of \cite{G16a} found for a statistics where the out- and in-regimes are not discriminated, as it should since, for $u\to1$, the 
requirement to be in the in-regime is not constraining ($p^{\rm in}(1)=1$). 
The probability density for $X(k,k\, u)$ for increasing values of $u$ and its universal limit above when $u\to 1$ are displayed in figure \ref{fig:probainLrescaled}.
\vskip .3cm
It is interesting to measure the relative contribution of the out- and in-regimes to the ``total" probability density 
for the rescaled length $L(d)$, i.e.\ the probability density obtained irrespectively of whether ${\mathcal V}(d)$ is finite or not, namely
\begin{equation*}
D(L,u)\equiv \lim_{k\to \infty} \frac{1}{dL} P_{k}\left(L\leq L(k\, u) < L+dL\right)=
D^{\rm out}(L,u)+D^{\rm in}(L,u)\ .
\end{equation*}
\begin{figure}
\begin{center}
\includegraphics[width=13cm]{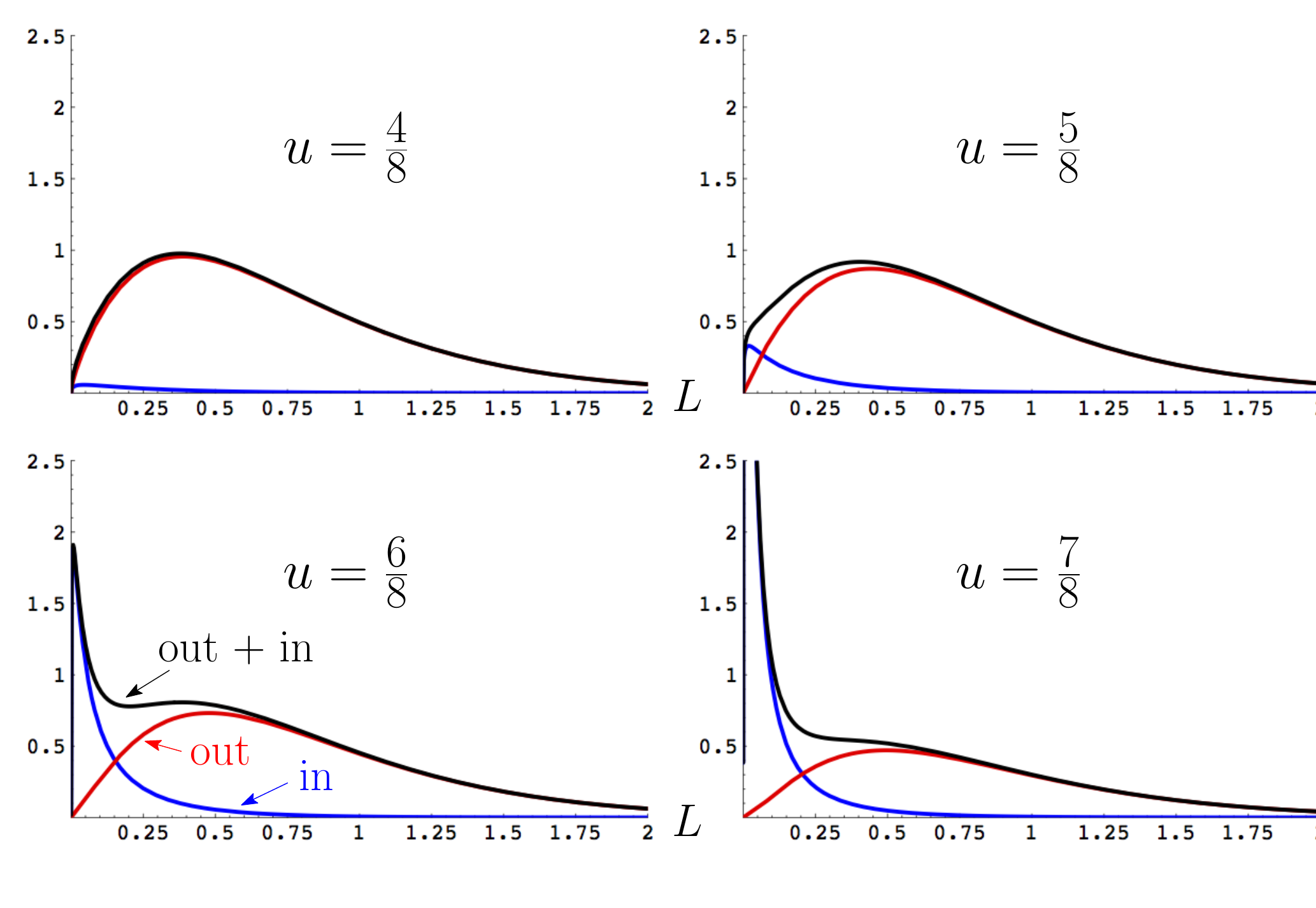}
\end{center}
\caption{
The relative contribution of the probability densities $D^{\rm out}(L,u)$ (in red) and $D^{\rm in}(L,u)$ (in blue) to the total probability density 
$D(L,u)=D^{\rm out}(L,u)+D^{\rm in}(L,u)$ (in black) as a function of $L$ for the indicated four values of $u$.
}
\label{fig:pLtot}
\end{figure}
The reader will easily check that the expression for $D(L,u)$ resulting from the explicit forms \eqref{eq:probainoutL} matches precisely the expression
for $D(L,u)$ found in \cite{G16a}, as it should. We have represented in figure  \ref{fig:pLtot} the probability density $D(L,u)$ for various values of $u$ as well as its two components
$D^{\rm out}(L,u)$ and $D^{\rm in}(L,u)$. As expected, $D(L,u)$ is dominated by the contribution of the out-regime at small enough $u$ (in practice up to $u\sim 1/2$) and 
by starts feeling the in-regime contribution when $u$ approaches $1$. This latter contribution moreover dominates the $u\to 1$ limit for small $L$. In particular, the appearance 
in $D(L,u)$ of a peak around $L=0$ when $u$ is large enough, 
which was observed in \cite{G16a} but remained quite mysterious is simply explained by the domination of the in-regime for $u\to 1$. No such peak
ever appears in the contribution $D^{\rm out}(L,u)$ of the out-regime.

From the laws \eqref{eq:probainoutL}, we may also compare the expectation value of $L(d)$ in the out- and in-regime to that obtained whithout conditioning: we
have respectively
\begin{equation}
\begin{split}
&\hskip -.5cm \lim_{k\to \infty} E_{k}\left(L(k\, u)\Big|{\mathcal V}(k\, u)\ \hbox{finite}\right)
=\frac{\lim
\limits_{k\to \infty}E_{k,k\, u}^{\rm out}\big(L(k\, u)\big)}{p^{\rm out}(u)}=\frac{3 c \left(4\!+\!4u\!-\!21 u^6\!+\!17 u^7\!-\!4 u^8\right)}{2 \left(4\!-\!7 u^6\!+\!3 u^7\right)}\ ,\\
&\hskip -.5cm\lim_{k\to \infty} E_{k}\left(L(k\, u)\Big|{\mathcal V}(k\, u)\ \hbox{infinite}\right)
=\frac{\lim\limits_{k\to \infty}E_{k,k\, u}^{\rm in}\big(L(k\, u)\big)}{p^{\rm in}(u)}=\frac{3 c (9\!-\!4 u) (1\!-\!u)}{2 (7\!-\!3 u)}\ , \\
&\hskip -.5cm\lim_{k\to \infty} E_{k}\left(L(k\, u)\right)
=p^{\rm out}(u)\ \frac{3 c \left(4+4u-21 u^6+17 u^7-4 u^8\right)}{2 \left(4-7 u^6+3 u^7\right)}+p^{\rm in}(u)\ \frac{3 c (9-4 u) (1-u)}{2 (7-3 u)}\\
& \hskip 1.9cm =\frac{3}{2} c \left(1+u-3u^6+u^7\right)\ ,\\
\end{split}
\label{eq:expectperim}
\end{equation}
where the last expression matches the result of \cite{G16a}.

\begin{figure}
\begin{center}
\includegraphics[width=12cm]{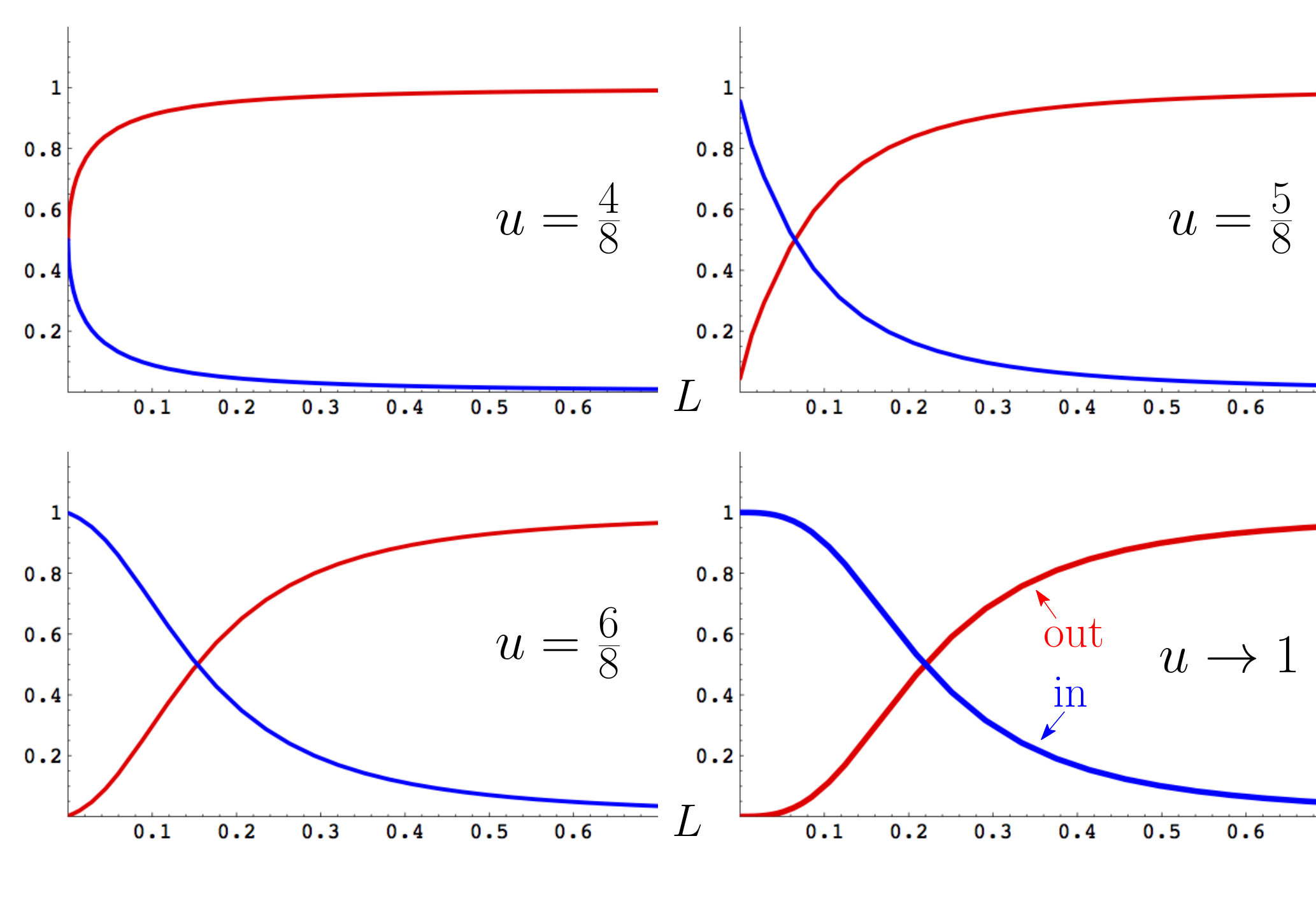}
\end{center}
\caption{The probabilities $\pi^{\rm out}(L,u)$ (in red) and $\pi^{\rm in}(L,u)$ (in blue) to be in the out- or in the in-regime,
knowing the value $L$ of the rescaled perimeter $L(k\, u)$ for fixed $u$ and in the limit $k\to \infty$. This probabilities
are represented as a function of $L$ for the indicated values of $u$.}
\label{fig:outorin}
\end{figure}

\begin{figure}
\begin{center}
\includegraphics[width=12cm]{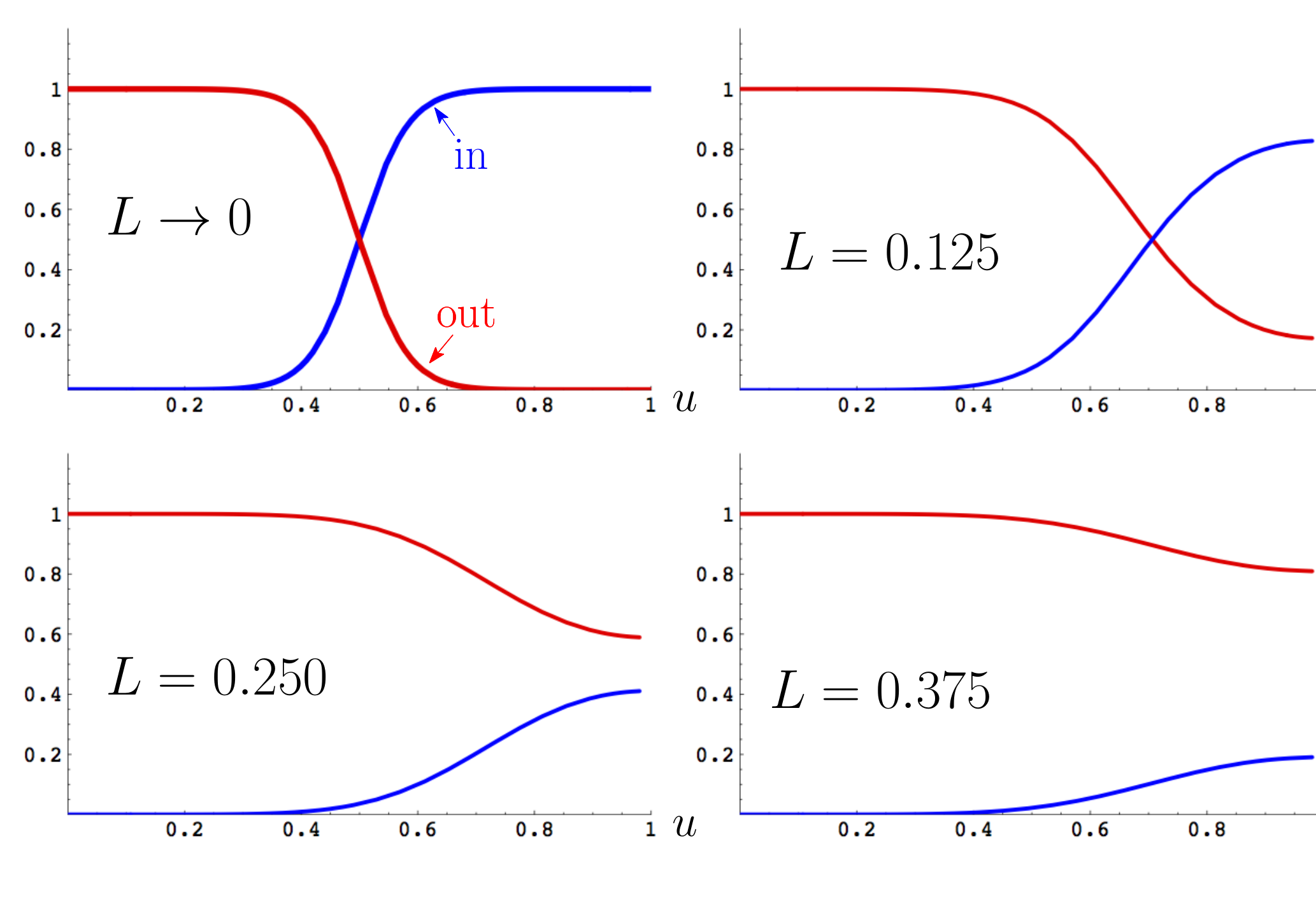}
\end{center}
\caption{The probabilities $\pi^{\rm out}(L,u)$ (in red) and $\pi^{\rm in}(L,u)$ (in blue) to be in the out- or in the in-regime,
knowing the value $L$ of the rescaled perimeter $L(k\, u)$ for fixed $u$ and in the limit $k\to \infty$. This probabilities
are represented as a function of $u$ for the indicated values of $L$.}
\label{fig:Uoutorin}
\end{figure}

To end this section, let us discuss the probability $\pi^{\rm out}(L,u)$ (resp. $\pi^{\rm in}(L,u)$ to be in the out-regime 
(resp.\ in the in-regime), knowing that the rescaled length $L(d)$ is equal to $L$ (with as before $u=d/k fixed$), namely
\begin{equation*}
\pi^{\rm out}(L,u)=\lim_{k\to \infty} P_{k}\left({\mathcal V}(k\, u)\ \hbox{finite}\Big| L(k\, u)=L\right)
= \frac{D^{\rm out}(L,u)}{D(L,u)}=1-\pi^{\rm in}(L,u)\ .
\end{equation*}
We have plotted in figure~\ref{fig:outorin} the quantities $\pi^{\rm out}(L,u)$ and $\pi^{\rm in}(L,u)$ as a function of $L$ for
various values of $u$. For $u\to 0$, we have $\pi^{\rm out}(L,0)=1$ and  $\pi^{\rm in}(L,0)=0$ irrespectively of $L$.
For $u\to 1$, we have the limiting expression: 
\begin{equation*}
\pi^{\rm out}(L,1)=\frac{28 L^3}{6 c^3+3 L c^2+28 L^3}=1-\pi^{\rm in}(L,1)\ .
\end{equation*}
Figure~\ref{fig:Uoutorin} displays the same probabilities $\pi^{\rm out}(L,u)$ and $\pi^{\rm in}(L,u)$, now as a function of $u$ for
various values of $L$. For $L\to 0$, we have the limiting expression
\begin{equation*}
\pi^{\rm out}(0,u)=\frac{(1-u)^6}{\left(1-2u+2 u^2\right) \left(1-4u+5 u^2-2u^3+u^4\right)}=1-\pi^{\rm in}(0,u) 
\end{equation*}
(note the remarkable symmetry $\pi^{\rm out}(0,u)=\pi^{\rm in}(0,1-u)$).
For $L\to \infty$, $\pi^{\rm out}(L,u)$ tends to $1$  and  $\pi^{\rm in}(L,u)$ to $0$, irrespectively of $u$.

\subsection{Joint law for the rescaled perimeter and volume in the out-regime}
\label{sec:jointlaw}
Our third result concerns the joint law for the hull perimeter ${\mathcal L}(d)$ and the hull volume ${\mathcal V}(d)$. Of course
such law is non-trivial only if the hull volume is finite, i.e.\ if we condition the configurations to be in the out-regime. From now on, all our results
will thus be \emph{conditioned to be in the out-regime}. For Large $d$,  ${\mathcal V}(d)$ scale as $d^4$ and we therefore introduce the rescaled volume
\begin{equation*}
V(d)\equiv \frac{{\mathcal V}(d)}{d^4}\ .
\end{equation*}
Our main result is the following expectation value
\begin{equation}
\begin{split}
&
 \lim_{k\to \infty} E_{k}\left(e^{-\sigma\, V(k\, u)-\tau\, L(k\, u)}\Big|{\mathcal V}(k\, u)\ \hbox{finite}\right)
=\frac{(1-u)^6}{u^3\, p^{\rm out}(u)} \times \frac{(f \sigma )^{3/4} \cosh \left(\frac{1}{2}(f \sigma)^{1/4}\right)}{8 \sinh^3 \left(\frac{1}{2}(f \sigma)^{1/4}\right)}
\\ & \hskip 10.cm \times M\big(\mu(\sigma,\tau,u)\big) \ ,\\
&\\
& \hbox{where}\ M(\mu)= \frac{1}{\mu ^4}\left(3 \mu ^2-5 \mu +6+\frac{4 \mu ^5+16 \mu ^4-7 \mu ^2-40 \mu -24}{4 (1+\mu )^{5/2}}\right)\\
&\hbox{and}\ \mu(\sigma,\tau,u)=\frac{(1-u)^2}{u^2}\, \left(
c\, \tau\,  +\frac{\sqrt{f \sigma }}{4}\,  \left(\coth^2 \left(\frac{1}{2}(f \sigma)^{1/4}\right)-\frac{2}{3}\right)\right)-1\ ,\\
\end{split}
\label{eq:expsigmatau}
\end{equation}
with $p^{\rm out}(u)$ as in \eqref{eq:poutin} and where $f$ is a normalization factor given in cases (i), (ii) and (iii) respectively by
\begin{equation}
\hbox{(i):}\  f=36\  , \qquad \hbox{(ii):}\  f=192\  , \qquad \hbox{(iii):}\  f=16\ .
\label{eq:fval}
\end{equation}
Setting $\tau=0$ and expanding at first order in $\sigma$, we immediately deduce that, in particular
\begin{equation*}
\hskip -.8cm \lim_{k\to \infty}  E_{k}\left(V(k\, u)\Big|{\mathcal V}(k\, u)\ \hbox{finite}\right)=\frac{f }{480}\frac{\left(20+12 u-77 u^6+57 u^7-12 u^8\right)}{ \left(4-7 u^6+3 u^7\right)}\ ,
\end{equation*}
a quantity which increases from $f/96$ at $u=0$ to $7f/480$ at $u=1$.
\vskip .3cm
The expectation value \eqref{eq:expsigmatau} above has a simple limit when $u\to 0$, namely
\begin{equation*}
\begin{split}
& \hskip -1.2cm \lim_{u\to 0}\Big( \lim_{k\to \infty} E_{k}\left(e^{-\sigma\, V(k\, u)-\tau\, L(k\, u)}\right)\Big)\\
& \hskip 3.cm =\frac{(f \sigma)^{3 /4}\cosh\left(\frac{1}{2}(f \sigma)^{1/4}\right)}{8 \sinh^3\left(\frac{1}{2}(f \sigma)^{1/4}\right)\left(c\, \tau\,  +\frac{\sqrt{f \sigma }}{4}\,  \left(\coth^2 \left(\frac{1}{2}(f \sigma)^{1/4}\right)-\frac{2}{3}\right)\right)^{3/2}}\\
\end{split}
\end{equation*}
(note that the condition that ${\mathcal V}(d)$ is finite is automatically satisfied in the limit $u\to 0$ since $p^{\rm out}(0)=1$). For $\tau=0$, this expression simplifies into
\begin{equation*}
\lim_{u\to 0}\Big( \lim_{k\to \infty} E_{k}\left(e^{-\sigma\, V(k\, u)}\right)\Big)=\frac{\cosh\left(\frac{1}{2}(f \sigma)^{1/4}\right)}{ \sinh^3\left(\frac{1}{2}(f \sigma)^{1/4}\right)\left(\coth^2 \left(\frac{1}{2}(f \sigma)^{1/4}\right)-\frac{2}{3}\right)^{3/2}}
\end{equation*}
and we recover here a result by Curien and Le Gall \cite{CLG14b}\footnote{The expression ${3}^{3/2} \cosh \left( (2\sigma)^{1/4}\, s/\sqrt{8/3}\right) \left(\cosh ^2\left((2\sigma)^{1/4}\, s/\sqrt{8/3}\right)+2\right)^{-3/2}$ of \cite{CLG14b} is indeed fully equivalent
under the correspondence $s=(2f)^{1/4}/\sqrt{3}$.}, in agreement with the equivalence principle \eqref{eq:utozero}. When $u\to 1$, we get another interesting limit
\begin{equation*}
\begin{split}
& \hskip -1.2cm \lim_{k\to \infty} E_{k}\left(e^{-\sigma\, V(k)-\tau\, L(k)}\Big|{\mathcal V}(k)\ \hbox{finite}\right)\\
& \hskip 3.cm
=\frac{(f \sigma)^{3 /4}\cosh\left(\frac{1}{2}(f \sigma)^{1/4}\right)}{8 \sinh^3\left(\frac{1}{2}(f \sigma)^{1/4}\right)\left(c\, \tau\,  +\frac{\sqrt{f \sigma }}{4}\,  \left(\coth^2 \left(\frac{1}{2}(f \sigma)^{1/4}\right)-\frac{2}{3}\right)\right)^{5/2}}\ .\\
\end{split}
\end {equation*}
\vskip .3cm
Performing an inverse Laplace transform on the variable $\tau$, we may extract from \eqref{eq:expsigmatau} the expectation value of $e^{-\sigma\, V(d)}$ knowing
the value $L$ of $L(d)$ in the out-regime. We find (see Section \ref{sec:expljointlaw} for details) that
\begin{equation}
\begin{split}
&\hskip -1.2cm \lim_{k\to \infty} E_{k}\left(e^{-\sigma\, V(k\, u)}\Big|{\mathcal V}(k\, u)\ \hbox{finite and}\ L(k\, u)
 =L\right) \\
& =\frac{1}{8} e^{-\frac{L}{c} \left(\frac{\sqrt{f \sigma }}{4} \left(\coth ^2\left(\frac{1}{2}(f \sigma)^{1/4}\right)-\frac{2}{3}\right)-1\right)} (f \sigma )^{3/4}
\frac{   \cosh \left(\frac{1}{2}(f \sigma )^{1/4}\right)}{ \sinh^3\left(\frac{1}{2}(f \sigma )^{1/4}\right)}\ .\\
   \end{split}
   \label{eq:VcondL}
\end{equation}
Note that this quantity turns out to be \emph{independent of $u$} and is thus equal to its limit for $u\to 0$. In agreement with
the equivalence \eqref{eq:utozero}, our result thus reproduces, \emph{now for any $u$},  the expression found by 
M\'enard in Ref.~\cite{Men16} in a limit where $k\to \infty$ before $d$ becomes large. 

We have in particular
\begin{equation*}
\lim_{k\to \infty} E_{k}\left(V(k\, u)\Big|{\mathcal V}(k\, u)\ \hbox{finite and}\ L(k\, u)=L\right)
 =\frac{f (c+L)}{240 c}
\end{equation*}
independently of $u$. The fact that the law for the rescaled volume $V(d)$, knowing the rescaled perimeter $L(d)$, is independent of $u$ 
is not so surprising. Indeed, $u$ measures the distance $k=d/u$ from the origin at which the marked vertex $x_1$ lies. Once the perimeter ${\mathcal L}(d)$
is fixed, the hull, whenever finite, is insensitive to the position of the second marked vertex. The law for its volume ${\mathcal V}(d)$ depends only on $d$ and ${\mathcal L}(d)$,
and, by simple scaling, it translates into a law for the rescaled volume $V(d)$ depending on the rescaled perimeter $L(d)$ only.
Note that, on the other hand, fixing the hull perimeter ${\mathcal L}(d)$ has some influence on the possible choices for the position 
of $x_1$ as a function of its distance $k=d/u$ from the origin $x_0$. This in return explains why the law for $L(d)$ and consequently that for $V(d)$ in the out-regime both depend on
$u$ for fixed $k$, as displayed in \eqref{eq:expsigmatau}.
 
\section{Derivation of the results: the strategy}
\label{sec:strategy}
Let us now come to the derivation of our results and explain the strategy behind our calculations. To simplify the discussion, we will focus here on the family (i) of quadrangulations.
The cases (ii) of triangulations and that (iii) of Eulerian triangulations are amenable to exactly the same type of treatment and we will briefly discuss
them in Section~\ref{sec:Other} below.

\subsection{Generating functions}
\label{sec:genfunc}
The main ingredient is the generating function $G(k,d,g,h,\alpha)$ of planar $k$-pointed-rooted quadrangulations, enumerated with a weight 
\begin{equation*}
g^{N-{\mathcal V}(d)}\ h^{{\mathcal V}(d)}\ \alpha^{{\mathcal L}(d)}\ ,
\end{equation*}
where $N$ is the total number of faces, and ${\mathcal L}(d)$ and ${\mathcal V}(d)$ are respectively the perimeter and volume of the hull at distance $d$
(note that ${\mathcal V}(d)\leq N$ by definition and we assume $k\geq 3$ and $2\leq d\leq k-1$).
To define precisely the hull at distance $d$, we use the construction discussed in \cite{G16a}. Then $G(k,d,g,h,\alpha)$ may be given an explicit expression as follows:
we use for the weights $g$ and $h$ the parametrization
\begin{equation}
g=\frac{x(1+x+x^2)}{(1+4x+x^2)^2}\  , \qquad h=\frac{y(1+y+y^2)}{(1+4y+y^2)^2}\ ,
\label{eq:gh}
\end{equation}
with $x$ and $y$ real between $0$ and $1$ (so that the generating function is well-defined for real $g$ and $h$ in the range $0\leq g,h \leq 1/12$).
We also introduce the quantity
\begin{equation*}
T_\infty(z)= \frac{z(1+4z+z^2)}{(1+z+z^2)^2}\ ,
\end{equation*}
where $z$ will be taken equal to $x$ or $y$ depending on the formula at hand.
We have, from \cite{G16a},
\begin{equation}
G(k,d,g,h,\alpha)= \underbrace{ \mathcal{K}\big( \mathcal{K}\big(\cdots \big( \mathcal{K}\big(}_{k-d\ \hbox{\scriptsize times}}\alpha^2\, T_{d}(y)\big)\big)\big)\big)- \underbrace{ \mathcal{K}\big( \mathcal{K}\big(\cdots \big( \mathcal{K}\big(}_{k-d\ \hbox{\scriptsize times}}\alpha^2\, T_{d-1}(y)\big)\big)\big)\big)\ ,
\label{eq:Gform}
\end{equation}
where $T_d(y)$ is defined by
\begin{equation*}
T_d(y)= T_\infty(y)\, \frac{(1-y^{d-1})(1-y^{d+4})}{(1-y^{d+1})(1-y^{d+2})}\ ,
\end{equation*}
and ${\mathcal K}\equiv \mathcal{K}(x)$ is an operator (depending on $x$ only), which satisfies the relation (which fully determines it):
\begin{equation}
\mathcal{K}\left(
T_\infty(x)\,  \frac{(1-\lambda\, x^{-1})(1-\lambda\, x^{4})}{(1-\lambda\, x)(1-\lambda\,  x^{2})}\right)=T_\infty(x)\, \frac{(1-\lambda)(1-\lambda\, x^{5})}{(1-\lambda\, x^{2})(1-\lambda\,  x^{3})}
\label{eq:propK}
\end{equation} 
for any arbitrary\footnote{In practice, as explained in \cite{G16a}, $\lambda$ must be small enough and this condition precisely dictates the branch 
of solution chosen in \eqref{eq:lambdaexpr} below.} $\lambda$ .
\begin{figure}
\begin{center}
\includegraphics[width=12cm]{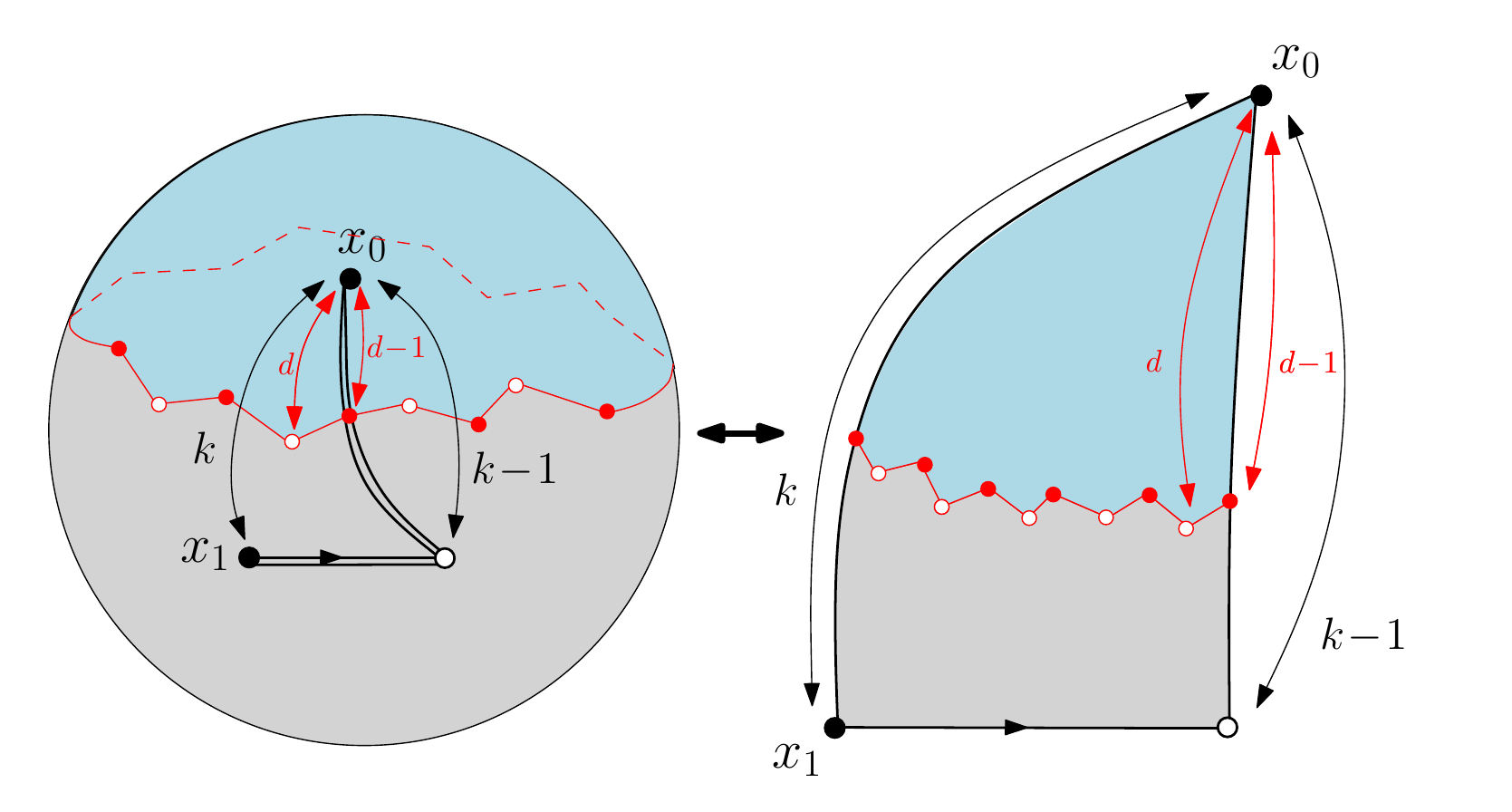}
\end{center}
\caption{A schematic picture of the bijection between a $k$-pointed-rooted planar quadrangulation (left) and a $k$-slice (right), as obtained by cutting
the quadrangulation along the leftmost shortest path from $x_1$ to $x_0$ (taking the root-edge of the map as first step). The light blue and light gray domains
are supposedly filled with faces of degree four. Left: the separating line at distance
$d$ (i.e.\ visiting alternately vertices at distance $d$ and $d-1$ from $x_0$) delimits the hull at distance $d$ (top part in light blue). Right: 
the image of this line connects the right- and left-boundaries of the $k$-slice and delimits an upper part containing $x_0$ (in light blue), which is the image of the hull,
from a lower part containing $x_1$.}
\label{fig:dividing}
\end{figure}
\begin{figure}
\begin{center}
\includegraphics[width=12cm]{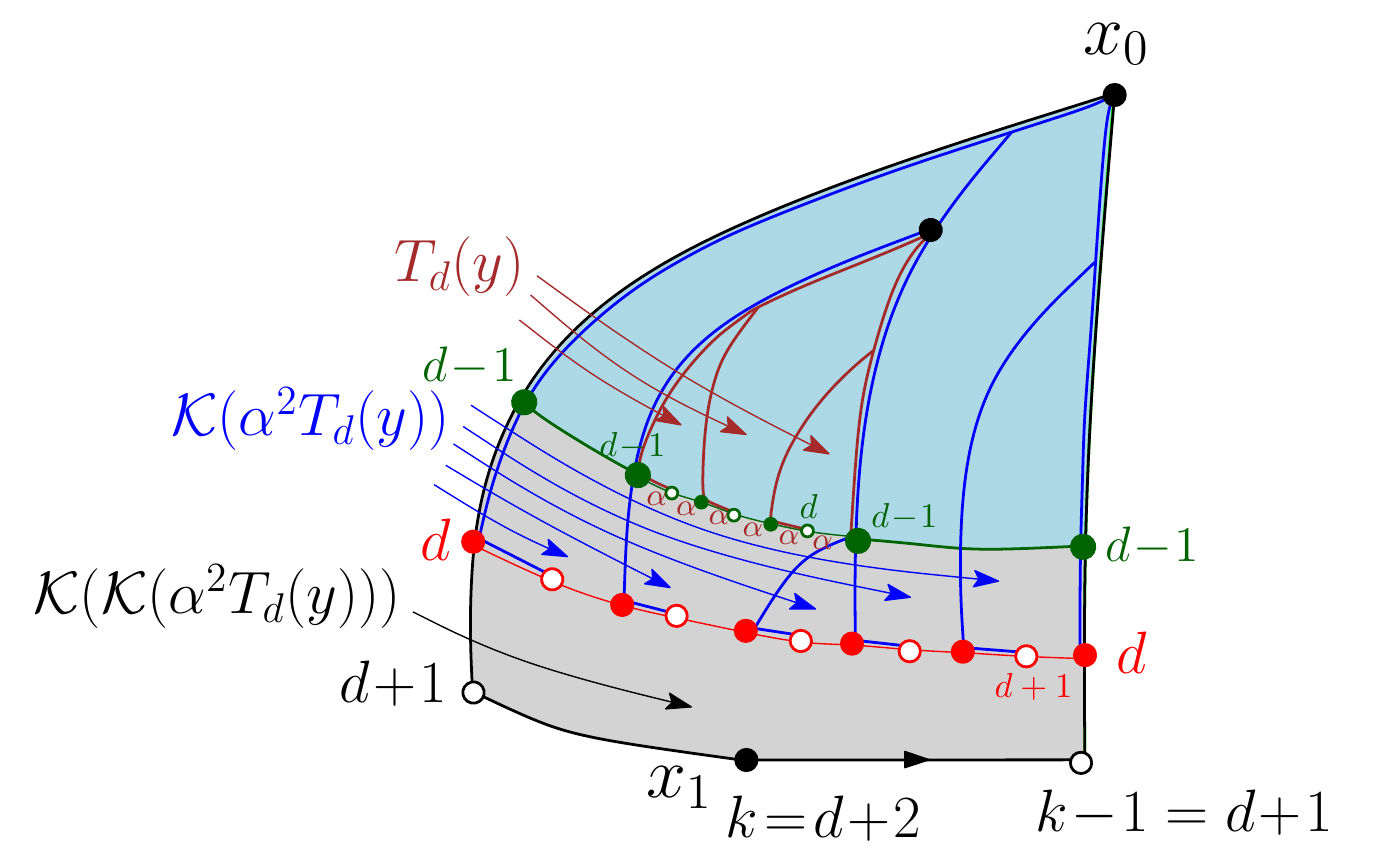}
\end{center}
\caption{A schematic picture of the successive decompositions of a $k$-slice obtained by cutting along $(k-d)$ successive separating lines at respective distance 
$d$, $d+1$, $\cdots$, $k-1$ from $x_0$ and, for the $m$-th such line ($1\le m\leq k-d$), by cutting along the leftmost 
shortest paths to $x_0$ from the $\mathcal{L}(d+m-1)/2$ vertices on this line lying at distance $d+m-1$ from $x_0$ (see text). Here $k=d+2$ and for the level $m=1$,
we represented only the leftmost shortest paths (in brown) lying within one particular sub-slice delimited by the leftmost shortest paths (in blue) at level $2$.}
\label{fig:constrhull}
\end{figure}

The origin of the above formula \eqref{eq:Gform} can be found in Refs.~\cite{G15b} and \cite{G16a}. We invite the reader to consult these references 
for details. Let us still briefly discuss the underlying decomposition 
of $k$-pointed-rooted quadrangulations on which the formula is based. As displayed in figure~\ref{fig:dividing}, a $k$-pointed-rooted quadrangulation
may be unwrapped into what is called a \emph{$k$-slice} by cutting it along some particular path of length $k$, namely the leftmost among shortest paths (along edges of the map)
from $x_1$ to $x_0$ having the root edge (i.e.\ the edge joining $x_1$ to its chosen neighbor at distance $k\!-\!1$ from $x_0$) as first step. The resulting $k$-slice
has a a left- and a right-boundary of respective lengths $k$ and $k\!-\!1$ linking the image of the root-edge in the $k$-slice (the so-called slice \emph{base}) to the image 
of $x_0$ (the so-called slice \emph{apex}). The passage from the $k$-pointed-rooted quadrangulation to the $k$-slice is a bijection so $G(k,d,g,h,\alpha)$
may also be viewed as the generating function for $k$-slices with appropriate weights.
The hull boundary at distance $d$ on the quadrangulation becomes a simple dividing line which links the right-and left-boundaries of the $k$-slice and separates it 
into an upper part, corresponding to the hull at distance $d$ in the original map and a complementary lower part. The dividing line visits alternately vertices at distance $d\!-\!1$
and $d$ from the apex, starting at the unique vertex along the right-boundary at distance $d\!-\!1$ from the apex and ending at the unique vertex along the left-boundary 
at distance $d\!-\!1$ from the apex. 
The upper part can be decomposed into a number of $d'$-slices with $d'\leq d$ by cutting it along the leftmost shortest paths to $x_0$ starting
from all the vertices at distance $d$ from $x_0$ along the dividing line. Since these vertices represent half of the vertices along the dividing line, the number
of $d'$-slices is precisely $\mathcal{L}(d)/2$ (note that $\mathcal{L}(d)$ is necessarily even for quadrangulations). Each of these slices 
is enumerated by a quantity $T_d(y)$ equal to the \emph{generating function of $d'$-slices with $2\leq d'\leq d$ and a weight $h$ per face} (see \cite{G15b} for
a precise definition), where $h$ and $y$ are related via \eqref{eq:gh}. The expression for $T_d(y)$ is that given just above, as computed in \cite{G15b}.
The juxtaposition of the $d'$-slices results in a total weight $(T_d(y))^{\mathcal{L}(d)/2}$ but, in order to impose that
the maximum value of $d'$ for all the $d'$-slices is actually \emph{exactly equal to} $d$, we must eventually subtract the weight of those configurations where 
all $d'$ would be less than $d\!-\!1$, namely $(T_{d-1}(y))^{\mathcal{L}(d)/2}$. Incorporating the desired weight $\alpha^{\mathcal{L}(d)}$, the generating function
of the upper part eventually reads
\begin{equation}
\left(\alpha^2\, T_d(y)\right)^{\mathcal{L}(d)/2}-\left(\alpha^2 T_{d-1}(y)\right)^{\mathcal{L}(d)/2}
\label{eq:upperpart}
\end{equation}
for the contribution of those configurations having a fixed value $\mathcal{L}(d)$ of the hull perimeter at distance $d$. 
This explains why the expression of $G(k,d,g,h,\alpha)$ is a difference of two terms, corresponding to the action of 
the same operator $\mathcal{K}^{\circ^{(k-d)}}$ on $\alpha^2\, T_d(y)$ and $\alpha^2\, T_{d-1}(y)$ respectively. 

To understand the origin of this operator, 
which incorporates
the contribution of the lower part, we proceed by recursion upon drawing the images of the successive hull boundaries at distance $d\!+\!1$, $d\!+\!2, \cdots$ 
until we reach the hull boundary at distance $k$ which reduces to the line of length $\mathcal{L}(k)=2$ formed by the concatenation of root-edge of the $k$-slice and
the first edge (starting from $x_1$) of the left-boundary (see figure \ref{fig:constrhull} where $k=d+2$). 
Looking at the hull boundary at distance $d\!+\!1$, we perform the same decomposition of the part above this boundary as we did before, by splitting
it into $\mathcal{L}(d+1)/2$ slices upon cutting along
the leftmost shortest paths to the apex starting from all the boundary vertices at distance $d+1$ (blue lines in figure \ref{fig:constrhull}). 
This creates $d''$-slices ${\mathcal S}_i$, $i=1,\cdots \mathcal{L}(d+1)/2$, each of them satisfying $d''\leq d+1$ and encompassing a number $\mathcal {L}_i$ of the previous $d'$-slices.
These $\mathcal {L}_i$ $d'$-slices contribute a weight $\left(\alpha^2\, T_d(y)\right)^{\mathcal{L}_i/2}$
to the first term in \eqref{eq:upperpart} (with $\mathcal{L}(d)=\sum_{i=1}^{\mathcal{L}(d+1)/2} \mathcal{L}_i$) while the generating function for the part of the slice
 ${\mathcal S}_i$ lying below the hull boundary at distance $d$, which \emph{depends only on the (half-)length $\mathcal{L}_i/2$} may 
 be written has  $[T^{\mathcal{L}_i/2}] \mathcal {K}(T)$ for some operator 
 $\mathcal {K}(T)$ depending on $g$ only (or equivalently on $x$ via \eqref{eq:gh}). This operator was computed in \cite{G15b} and, as explained in \cite{G16a},  satisfies the 
 property \eqref{eq:propK} above. Summing over all values of $\mathcal{L}(d)$, hence on all values of $\mathcal{L}_i$, each of the $d''$-slices contributes
 a weight 
 \begin{equation*}
 \mathcal{K}\left(\alpha^2\, T_d(y)\right)
 \end{equation*} 
 to the sum over $\mathcal{L}(d)$ of the first term in \eqref{eq:upperpart}. Taking into account the $\mathcal{L}(d+1)/2$ $d''$-slices, 
 we end up with a contribution
\begin{equation*}
\left(\mathcal{K}\left(\alpha^2\, T_d(y)\right)\right)^{\mathcal{L}(d+1)/2}-\left(\mathcal{K}\left(\alpha^2 T_{d-1}(y)\right)\right)^{\mathcal{L}(d+1)/2}
\end{equation*}
for the part above the hull-boundary at distance $d+1$ of those configurations with a fixed value $\mathcal{L}(d+1)$ of the hull perimeter at distance $d+1$. 
Repeating the argument $k-d$ times immediately yields the desired expression \eqref{eq:Gform} since $\mathcal{L}(k)=2$ by construction.
\vskip .3cm
In order to have a more tractable expression, we may now perform explicitly the $k\!-\!d$ iterations of the operator ${\mathcal K}$ in \eqref{eq:Gform}. This leads 
immediately to
the more explicit formula
\begin{equation}
\begin{split}
& \hskip -1.cm G(k,d,g,h,\alpha)= H\big(k-d,x,\alpha^2T_d(y)\big)-H\big(k-d,x,\alpha^2T_{d-1}(y)\big)\\
& \hbox{where}\ H(k,x,T)=T_\infty(x)\,  \frac{\big(1-\lambda(x,T)\, x^{k-1}\big)\big(1-\lambda(x,T)\, x^{k+4}\big)}{\big(1-\lambda(x,T)\, x^{k+1}\big)\big(1-\lambda(x,T)\,  x^{k+2}\big)} \\
\end{split}
\label{eq:Gexpr}
\end{equation}
provided $\lambda(x,T)$ is defined through
\begin{equation*}
T_\infty(x)\,  \frac{\big(1-\lambda(x,T)\, x^{-1}\big)\big(1-\lambda(x,T)\, x^{4}\big)}{\big(1-\lambda(x,T) \, x\big)\big(1-\lambda(x,T) \,  x^{2}\big)}=T \ ,
\end{equation*}
namely
\vskip -1.cm
\begin{equation}
\begin{split}
&\\ &\\
&\hskip -1.cm \lambda(x,T)=\frac{T_\infty(x)(1\!+\!x^5)\!-x^2T (1+x)\!-\!\sqrt{\left(T_\infty(x)(1\!+\!x^5)\!-\!x^2T (1\!+\!x)\right)^2\!-\!4\, x^5(T_\infty(x)\!-\!T)^2}}
{2\, x^4\, (T_\infty(x)\!-\!T)}\ , \\
\end{split}
\label{eq:lambdaexpr}
\end{equation}
These latest expressions \eqref{eq:Gexpr} and \eqref{eq:lambdaexpr} will be our starting point for explicit calculations.
\vskip .3cm
A last quantity of interest is the generating function of $F(k,g)$ of planar $k$-pointed-rooted quadrangulations with a weight $g$ per face.
We have clearly $F(k,g)=G(k,d,g,g,1)$ for any $d\leq k-1$ and we easily obtain from the above formulas that $\lambda\big(x,T_d(x)\big)=x^d$ so that
\begin{equation*}
F(k,g)=T_\infty(x)\left( \frac{(1-x^{k-1})(1-x^{k+4})}{(1-x^{k+1})(1-x^{k+2})}-\frac{(1-x^{k-2})(1-x^{k+3})}{(1-x^{k})(1-x^{k+1})}\right)\ .
\end{equation*}

\subsection{Sending ${\boldsymbol N\to \infty}$: the out- and in-regimes}
\label{sec:extract}
Let us now explain how we can extract from the above generating functions results on the $N\to \infty$ limit, imposing that the
configurations are either in the out- or the in-regime. To simplify the notations, let us omit for a while the dependence of $G(k,d,g,h,\alpha)$
in $\alpha$, $k$ and $d$ and write $G(k,d,g,h,\alpha)=G(g,h)$, as well as $N-{\mathcal V}(d)=n_1$ and ${\mathcal V(d)}=n_2$.
We also denote by $G_{n_1,n_2}$ the coefficient $[g^{n_1}h^{n_2}]G(g,h)$. 
We are then interested in the large $N$ limit of the quantity
\begin{equation}
\sum_{n_1,n_2\atop n_1+n_2=N} G_{n_1,n_2}\ ,
\label{eq:sumn1n2}
\end{equation}
which we wish to extract from the knowledge of the generating function $G(g,h)$. As mentioned earlier, we also assume that when $N\to \infty$, the sum in \eqref{eq:sumn1n2} 
is dominated by two contributions, that with $n_1\to \infty$, $n_2$ staying finite,
which corresponds to what we called the out-regime, and that with $n_2\to \infty$, $n_1$ staying finite,
which we called the in-regime, while the contribution where both $n_1$ and $n_2$ become infinite simultaneously is algebraically suppressed for large $N$. 
To describe the out-regime, we must consider the $n_1\to \infty$ behavior of 
$G_{n_1,n_2}$ which is encoded in the singular behavior of $G(g,h)$ when $g$ reaches some critical value $g^*$ (the radius of convergence of the series in $g$,
possibly depending on $h$ and the other parameters).
Similarly, properties of the in-regime are encoded in the singular behavior of $G(g,h)$ when $h$ reaches some critical value $h^*$ (possibly depending on $g$ and the
other parameters).
For the generating function $G(k,d,g,h,\alpha)$ of interest, the singularities appear when either $x\to 1$ or $y\to 1$, irrespectively of $k$, $d$ and $\alpha$,
i.e., from \eqref{eq:gh}, for $g\to 1/12$ or $h\to 1/12$.
We therefore have $g^*=1/12$ (independently of $h$ and the other parameters) and $h^*=1/12$ (independently of $g$ and the other parameters). 
More precisely, we have expansions of the form
\begin{equation*}
\begin{split}
&G(g,h)=\mathfrak{g}_0(h)+\mathfrak{g}_2(h)(1-12 g)+\mathfrak{g}_3(h)(1-12 g)^{3/2}+O((1-12g)^2)\ ,\\
&G(g,h)=\tilde{\mathfrak{g}}_0(g)+\tilde{\mathfrak{g}}_2(g)(1-12 h)+\tilde{\mathfrak{g}}_3(g)(1-12 h)^{3/2}+O((1-12h)^2)\ .\\
\end{split}
\end{equation*}
(where all the functions implicitly depend on $k$, $d$ and $\alpha$). Note in particular that $G(1/12,h)=\mathfrak{g}_0(h)$ and $G(g,1/12)=\tilde{\mathfrak{g}}_0(g)$ 
are finite and that there are no square-root singularities.

Taking the term of order $h^{n_2}$ in the first expansion above, we deduce the singular part
\begin{equation*}
\left(\sum_{n_1} G_{n_1,n_2} g^{n_1}\right) \Big|_{\rm sing.}=[h^{n_2}] \mathfrak{g}_3(h) \times (1-12 g)^{3/2}
\end{equation*}
from which we deduce the large $n_1$ behavior
\begin{equation*}
G_{n_1,n_2} \underset{n_1\to \infty}{\sim}[h^{n_2}] \mathfrak{g}_3(h)\times \frac{3}{4}  \frac{12^{n_1}}{\sqrt{\pi} n_1^{5/2}}
\end{equation*}
so that 
\begin{equation*}
\sum_{n_2} G_{N-n_2,n_2} \underset{N \to \infty}{\sim}\frac{3}{4}  \frac{12^{N}}{\sqrt{\pi} N^{5/2}} \sum_{n_2} [h^{n_2}] \mathfrak{g}_3(h) \times 12^{-n_2}
= \frac{3}{4}  \frac{12^{N}}{\sqrt{\pi} N^{5/2}} \mathfrak{g}_3\left(\frac{1}{12}\right) \ .
\end{equation*}
This represents precisely \emph{the contribution of the out-regime} to the large $N$ limit of the sum \eqref{eq:sumn1n2}. If, more generally, we
wish to control the volume $\mathcal{V}(d)$ in the out-regime, we may consider
\begin{equation*}
\sum_{n_2} G_{N-n_2,n_2}\, \rho^{n_2} \underset{N \to \infty}{\sim}\frac{3}{4}  \frac{12^{N}}{\sqrt{\pi} N^{5/2}} \sum_{n_2} [h^{n_2}] \mathfrak{g}_3(h) \times \left(\frac{\rho}{12}\right)^{n_2}
= \frac{3}{4}  \frac{12^{N}}{\sqrt{\pi} N^{5/2}} \mathfrak{g}_3\left(\frac{\rho}{12}\right) \ .
\end{equation*}
By a similar argument, the in-regime contribution to the sum \eqref{eq:sumn1n2} behaves as 
\begin{equation*}
\sum_{n_1} G_{n_1,N-n_1} \underset{N \to \infty}{\sim} \frac{3}{4}  \frac{12^{N}}{\sqrt{\pi} N^{5/2}} \tilde{\mathfrak{g}}_3\left(\frac{1}{12}\right)\ .
\end{equation*}
For $F(g)\equiv F(k,g)$, we have an expansion of the form
\begin{equation*}
F(g)=\mathfrak{f}_0+\mathfrak{f}_2(1-12 g)+\mathfrak{f}_3(1-12 g)^{3/2}+O((1-12g)^2)
\end{equation*}
which yields the large $N$ estimate
\begin{equation*}
[g^N]F(g)  \underset{N \to \infty}{\sim} \frac{3}{4}  \frac{12^{N}}{\sqrt{\pi} N^{5/2}} \mathfrak{f}_3 \ .
\end{equation*}
By taking the appropriate ratios, we eventually deduce the large $N$ limit of the desired expectation values, namely
\begin{equation}
E_{k,d}^{\rm out}\left(\rho^{\mathcal{V}(d)}\, \alpha^{\mathcal{L}(d)}\right)
= \frac{\mathfrak{g}_3\left(\frac{\rho}{12},k,d,\alpha\right)}{\mathfrak{f}_3(k)}\  ,
\label{eq:asympout}
\end{equation}
where we re-introduced explicitly the dependence in $k$, $d$ and $\alpha$ of $\mathfrak{g}_3(h)\equiv \mathfrak{g}_3(h,k,d,\alpha)$ and $\mathfrak{f}_3\equiv \mathfrak{f}_3(k)$,
and 
\begin{equation}
E_{k,d}^{\rm in}\left(\alpha^{\mathcal{L}(d)}\right)
= \frac{\tilde{\mathfrak{g}}_3\left(\frac{1}{12},k,d,\alpha\right)}{\mathfrak{f}_3(k)}
\label{eq:asympin}
\end{equation}
(with the more explicit dependence $\tilde{\mathfrak{g}}_3(g)\equiv \tilde{\mathfrak{g}}_3(g,k,d,\alpha)$). 
This reduces our problem to estimating the quantities $\mathfrak{g}_3(h)$, $\tilde{\mathfrak{g}}_3(g)$ and $\mathfrak{f}_3$
from our explicit expressions for $G(k,d,g,h,\alpha)$ and $F(k,g)$.

\section{Derivation of the results: explicit calculations}
\label{sec:explicit}
Since the expression for $G(k,d,g,h,\alpha)$ is quite involved, explicit calculations may be difficult to perform in all generalities for 
finite $k$ and $d$ and some of our results will hold only in the limit of large $k$ and $d$. Still, the simplest questions may be solved exactly
for finite $k$ and $d$: this is the case for the probability to be in the out- or the in-regime, as we discuss now.

\subsection{Results at finite ${\boldsymbol k}$ and ${\boldsymbol d}$: the probability to be in the out- or in-regime}
\label{sec:explprobs}
If we wish to compute the probability to be in the out- or in-regime, we may set $\alpha=1$ and $\rho=1$ in \eqref{eq:asympout} and \eqref{eq:asympin}
since we do not measure 
the hull perimeter nor the hull volume. 
More precisely, we have
\begin{equation*}
\begin{split}
&P_k(\mathcal{V}(d)\ \hbox{finite})=E_{k,d}^{\rm out}(1)= \frac{\mathfrak{g}_3\left(\frac{1}{12},k,d,1\right)}{\mathfrak{f}_3(k)}\  ,\\
&P_k(\mathcal{V}(d)\ \hbox{infinite})=E_{k,d}^{\rm in}(1)= \frac{\tilde{\mathfrak{g}}_3\left(\frac{1}{12},k,d,1\right)}{\mathfrak{f}_3(k)}\  .\\
\end{split}
\end{equation*}
To compute $\mathfrak{f}_3(k)$, we set
\begin{equation*}
g=\frac{1}{12}(1-\epsilon^4)\leftrightarrow \epsilon=(1-12 g)^{1/4}\ ,
\end{equation*}
and, from the corresponding small $\epsilon$ expansion of $x$,
\begin{equation*}
 x=1-\sqrt{6}\, \epsilon +3\, \epsilon ^2-\frac{5}{2} \sqrt{\frac{3}{2}}\, \epsilon ^3+3\, \epsilon ^4-\frac{39}{16} \sqrt{\frac{3}{2}}\, \epsilon ^5+3\, \epsilon ^6-\frac{157}{64}
  \sqrt{\frac{3}{2}} \, \epsilon ^7+3 \epsilon ^8-\cdots\ ,
\end{equation*}
we immediately get the small $\epsilon$ expansion of $F(k,g)$. 
As expected, we find terms of order $\epsilon^0$, $\epsilon^4=(1-12 g)$ and $\epsilon^6=(1-12 g)^{3/2}$ but no term of
odd order in $\epsilon$ (as a consequence of the $x\to 1/x$ symmetry of all the formulas) and, more importantly, no term of order $\epsilon^2=(1-12 g)^{1/2}$.
The coefficient of $\epsilon^6$ in the expansion yields:
\begin{equation}
\mathfrak{f}_3(k)= \frac{4 \left(k^2+2 k-1\right) \left(5 k^4+20 k^3+27 k^2+14 k+4\right)}{35 k (k+1) (k+2)}\ .
\label{eq:f3}
\end{equation}
To compute $\mathfrak{g}_3\left(1/12,k,d,1\right)$, we have to consider the expansion of  $G(k,d,g,1/12,1)$ when $g\to 1/12$.
Note that setting $h=1/12$ amounts to setting $y=1$, in which case $T_d(y)$ simplifies into
\begin{equation*}
T_d(1)=\frac{2}{3}\ \frac{(d-1)(d+4)}{(d+1)(d+2)}\ .
\end{equation*}
We may easily compute the singularity of the function $H(k,x,T)$ (as defined in \eqref{eq:Gexpr} and \eqref{eq:lambdaexpr}) for $\epsilon\to 0$. 
Again, the leading singularity corresponds to the $\epsilon^6$ term and we find explicitly
\begin{equation}
\begin{split}
&
\hskip -1.cm H(k,x,T)\Big|_{\rm sing.}=\mathfrak{h}_3\big(k,Y(T)\big) (1-12 g)^{3/2}\\
&\hskip -1.cm \hbox{with}\ \mathfrak{h}_3(k,Y)= \frac{k}{840\,  Y \left((2 k+Y)^2-1\right)^2}\Big(
105 (k+Y)^8+420 \left(k^2-3\right) (k+Y)^6\\& -210 \left(k^4+6 k^2+49\right) (k+Y)^4-4 \left(75 k^6-567 k^4-1715 k^2-2273\right) (k+Y)^2\\
& -(k-5) (k-1) (k+1)
   (k+5) \left(15 k^4+138 k^2-217\right)\Big)\\
   & \hskip -1.cm \hbox{and}\ Y(T)=\sqrt{\frac{3 T-50}{3 T-2}}\ .
   \\
   \end{split}
   \label{eq:Hsing}
\end{equation}
Note that we have the particularly simple expression
\begin{equation*}
Y\big(T_d(1)\big)=2d+3
\end{equation*}
so that, from \eqref{eq:Gexpr}, 
\begin{equation*}
\mathfrak{g}_3\left(\frac{1}{12},k,d,1\right)=\mathfrak{h}_3(k-d,2d+3)-\mathfrak{h}_3(k-d,2d+1)
\end{equation*}
and
\begin{equation*}
\begin{split}
&\hskip -1.cm P_k({\mathcal V}(d)\ \hbox{finite})=\frac{1}{\mathfrak{f}_3(k)}\left(\mathfrak{h}_3(k-d,2d+3)-\mathfrak{h}_3(k-d,2d+1)\right)\\
& \hskip -1.cm = \frac{1}{\mathfrak{f}_3(k)}\Bigg(
\frac{1}{105 (2d+3) (k+1)^2 (k+2)^2}\times \\ & \hskip .3cm\Big((2 d+3) (k-1) (k+1) (k+2) (k+4) \left(15 k^4+90 k^3+237 k^2+306 k+140\right)\\ &\hskip .3cm - (2 k+3) (d-1) (d+1) (d+2) (d+4) \left(15 d^4+90 d^3+237 d^2+306 d+140\right)\Big)\\
& -\frac{1}{105 (2 d+1) k^2 (k+1)^2}\times \\ & 
\hskip .3cm\big((2 d+1) (k-2) k (k+1) (k+3) \left(15 k^4+30 k^3+57 k^2+42 k-4\right)\\ &\hskip .3cm  -(2k+1)(d-2) d (d+1) (d+3) \left(15 d^4+30 d^3+57 d^2+42 d-4\right)\big)
\Bigg)\\
\end{split}
\end{equation*}
with $\mathfrak{f}_3(k)$  as in \eqref{eq:f3} above.

\begin{figure}
\begin{center}
\includegraphics[width=7.5cm]{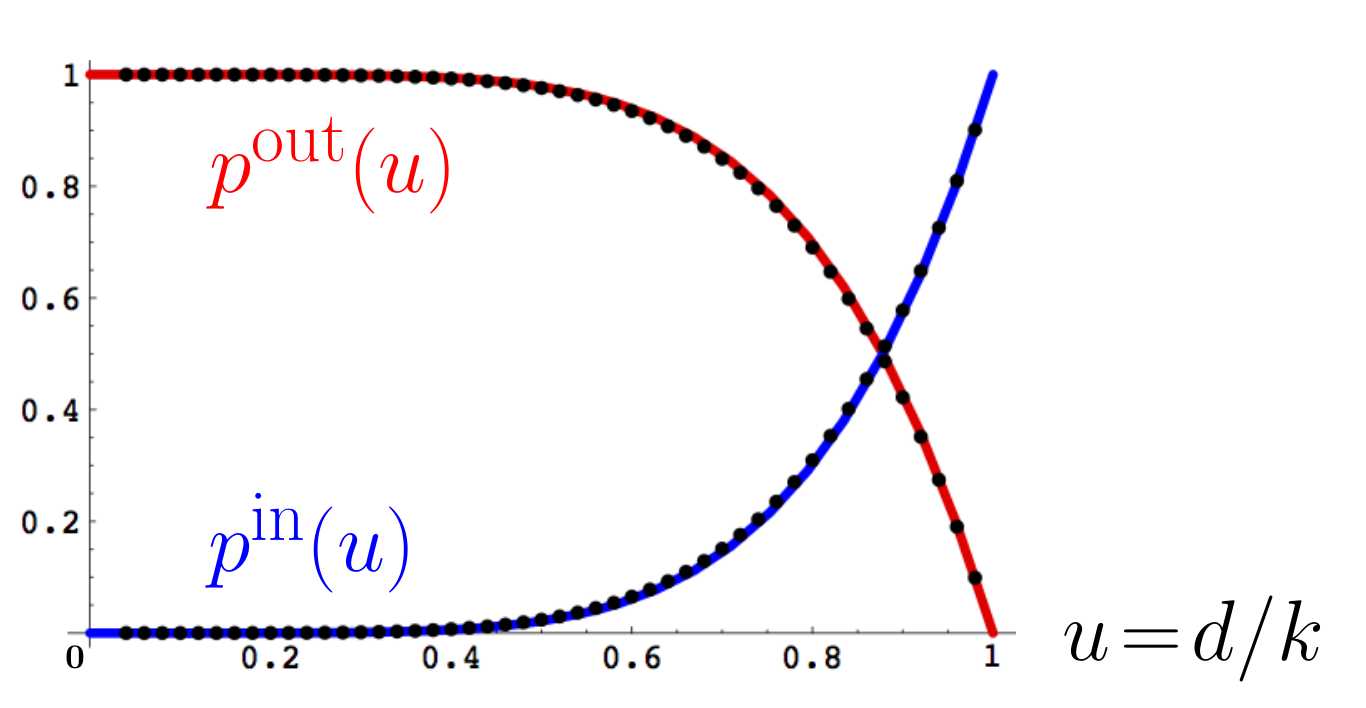}
\end{center}
\caption{A comparison between the probability $P_k({\mathcal V}(d)\ \hbox{finite})$ (respectively $P_k({\mathcal V}(d)\ \hbox{infinite})$) vs $d/k$ for 
$k=50$ and $2\leq d\leq 49$ and its limiting expression $p^{\rm out}(u)$ (respectively $p^{\rm in}(u)$) vs $u$, as given by \eqref{eq:poutin}.}
\label{fig:probasoutinfinite}
\end{figure}

Let us now compute the probability to be in the in-regime which requires the knowledge of $\tilde{\mathfrak{g}}_3\left(1/12,k,d,1\right)$. 
We thus have to consider the expansion of $G(k,d,1/12,h,1)$ when $h\to 1/12$.
Note that setting $g=1/12$ now amounts to setting $x=1$\footnote{More precisely, we must take the limit $x\to 1^-$.}, in which case 
we find the simple expression
\begin{equation}
H(k,1,T)=\frac{2}{3}\frac{\left(2 k+Y\left(T\right)\right)^2-25}{\left(2 k+Y\left(T\right)\right)^2-1}
\label{eq:HkoneT}
\end{equation}
with $Y(T)$ as in \eqref{eq:Hsing}. To compute the desired singularity, we now set
\begin{equation*}
h=\frac{1}{12}(1-\eta^4)\leftrightarrow \eta=(1-12 h)^{1/4}\ ,
\end{equation*}
so that 
\begin{equation*}
 y=1-\sqrt{6}\, \eta +3\, \eta ^2-\frac{5}{2} \sqrt{\frac{3}{2}}\, \eta ^3+3\, \eta ^4-\frac{39}{16} \sqrt{\frac{3}{2}}\, \eta ^5+3\,\eta ^6-\frac{157}{64}
  \sqrt{\frac{3}{2}} \, \eta ^7+3 \eta ^8-\cdots\ .
\end{equation*}
We have the expansion
\begin{equation}
\begin{split}
&\hskip -1.cm Y\big(T_d(y)\big)=(2d+3)-\frac{(d-1) (d+1) (d+2) (d+4) \left(9 d^2+27 d+10\right)}{30 (2 d+3)}\eta^4\\
& +
\frac{(d-1) (d+1) (d+2) (d+4) \left(15 d^4+90 d^3+237 d^2+306 d+140\right)}{210 (2 d+3)}\eta^6+\cdots \\
\end{split}
\label{eq:expTd}
\end{equation}
which yields eventually
\begin{equation*}
\begin{split}
&
\hskip -1.cm H\big(k,1,T_d(y)\big)\Big|_{\rm sing.}=\tilde{\mathfrak{h}}_3(k,d) (1-12 h)^{3/2}\\
&\hskip -1.cm \hbox{with}\ \tilde{\mathfrak{h}}_3(k,d)= \frac{(d-1) (d+1) (d+2) (d+4) \left(15 d^4+90 d^3+237 d^2+306 d+140\right) (2 d+2 k+3)}{105 (2 d+3) (d+k+1)^2 (d+k+2)^2}\ .
   \\
   \end{split}
\end{equation*}
We end up with 
\begin{equation*}
\begin{split}
&\hskip -1.cm P_k({\mathcal V}(d)\ \hbox{infinite})=\frac{1}{\mathfrak{f}_3(k)}\left(\tilde{\mathfrak{h}}_3(k-d,d)-\tilde{\mathfrak{h}}_3(k-d,d-1)\right)\\
& \hskip -1.cm = \frac{1}{\mathfrak{f}_3(k)}\Bigg(
\frac{(2 k+3)(d-1) (d+1) (d+2) (d+4) \left(15 d^4+90 d^3+237 d^2+306 d+140\right)}{105 (2 d+3) (k+1)^2 (k+2)^2}\Big)\\
& \hskip 1.cm-\frac{(2 k+1)(d-2) d (d+1) (d+3) \left(15 d^4+30 d^3+57 d^2+42 d-4\right) }{105 (2 d+1) k^2 (k+1)^2}
\Bigg)\ .\\
\end{split}
\end{equation*}
It is easily verified from their explicit expressions that 
\begin{equation*}
P_k({\mathcal V}(d)\ \hbox{finite})+P_k({\mathcal V}(d)\ \hbox{infinite})=1
\end{equation*}
for any fixed $k$ and $d$, as expected. This corroborates the absence of some regime other than the out- and in-regimes and
justifies a posteriori our statement that the contribution of configurations where both the hull and its complementary would have infinite volumes is negligible at large $N$. 
For $k\to\infty$ and $d\to \infty$ with $u=d/k$ fixed, we immediately obtain
\begin{equation*}
\begin{split}
& \lim_{k\to \infty} P_k({\mathcal V}(k\, u)\ \hbox{finite})=\frac{1}{4} \left(4-7 u^6+3 u^7\right)\ ,\\
&\lim_{k\to \infty} P_k({\mathcal V}(k\, u)\ \hbox{infinite})=\frac{1}{4} (7-3 u)\, u^6\ ,\\
\end{split}
\end{equation*}
which is precisely the announced result \eqref{eq:poutin}. Figure~\ref{fig:probasoutinfinite} shows a comparison 
between the limiting expressions $p^{\rm out}(u)$ and $p^{\rm in}(u)$ vs $u$ (as given by \eqref{eq:poutin}) and the 
corresponding finite $k$ and $d$ expressions (as given above)
$P_k({\mathcal V}(d)\ \hbox{finite})$ and $P_k({\mathcal V}(d)\ \hbox{infinite})$ vs $d/k$ for $k=50$ and $2\leq d\leq 49$.
\vskip .3cm
Another quantity which may be easily computed for finite $k$ and $d$ is the expectation value of the perimeter in the out-regime, 
$E_{k}\left({\mathcal L}(d)\Big|{\mathcal V}(d)\ \hbox{finite}\right)$, as well as that in the in-regime, $E_{k}\left({\mathcal L}(d)\Big|{\mathcal V}(d)\ \hbox{infinite}\right)$. 
Details of the computation are discussed in Appendix A.

\subsection{Law for the perimeter at large ${\boldsymbol k}$ and ${\boldsymbol d}$ in the out- and in-regimes}
\label{sec:explperimlaw}
To describe the statistics of the perimeter in the out-regime, we have to look at $G(k,d,g,h,\alpha)$ for arbitrary $\alpha$
and for $h=1/12$ (i.e.\ $y=1$). We consider here the large $k$ and $d$ limit by setting $d=k\, u$ (with $0\leq u\leq 1$) 
and letting $k\to \infty$. In this limit, ${\mathcal L}(d)$ growths like $d^2=(k\, u)^2$ and the large $k$ statistics of the perimeter
is captured by setting
\begin{equation*}
\alpha=e^{-\frac{\tau}{(k\, u)^2}}
\end{equation*}
with $\tau$ finite. From \eqref{eq:Gexpr} and \eqref{eq:Hsing}, we need the large $k$ behavior of $\mathfrak{h}_3(k,Y(T))$ for 
$T=\alpha^2 T_{d}(1)$ and  $T=\alpha^2 T_{d-1}(1)$ which involves the associated expansions of $Y(T)$, namely
\begin{equation*}
\begin{split}
&Y\left(e^{-2 \frac{\tau}{(k\, u)^2}} T_{k\, u}(1)\right)= 2 \sqrt{\frac{3}{3+\tau}}\, k\, u+ 3 \left(\frac{3}{3+\tau}\right)^{3/2}+O\left(\frac{1}{k}\right)\ ,\\
&Y\left(e^{-2 \frac{\tau}{(k\, u)^2}} T_{k\, u-1}(1)\right)= 2 \sqrt{\frac{3}{3+\tau}}\, k\, u+ \left(\frac{3}{3+\tau}\right)^{3/2}+O\left(\frac{1}{k}\right)\ .\\
\end{split}
\end{equation*}
Using the expansion
\begin{equation*}
\mathfrak{f}_3=\frac{4}{7} k^3+O\left(k^2\right)\ , 
\end{equation*}
we obtain\footnote{Note that, for $\tau=0$, $\frac{(1-u)^6 }{u^3} M\big(\mu(0,u)\big)=\frac{1}{4} \left(4-7 u^6+3 u^7\right)=p^{\rm out}(u)$ as it should.}
\begin{equation}
\begin{split}
& \hskip -1.cm \lim_{k\to \infty} E_{k,k\, u}^{\rm out}\left(e^{-\tau\,  L(k\, u)}\right)\\
&=\lim_{k\to \infty}
\frac{\mathfrak{h}_3\left(k\!-\!k\, u,2 \sqrt{\frac{3}{3+\tau}}\, k\, u\!+\! 3 \left(\frac{3}{3+\tau}\right)^{3/2} \right)\!-\!\mathfrak{h}_3\left(k\!-\!k\, u,2 \sqrt{\frac{3}{3+\tau}}\, k\, u\!+\! \left(\frac{3}{3+\tau}\right)^{3/2}\right)}{\frac{4}{7} k^3}\\
&=\frac{(1-u)^6 }{u^3} M\big(\mu(\tau,u)\big) \\
&\hskip -1.cm  \hbox{with}\ M(\mu)= \frac{1}{\mu ^4}\left(3 \mu ^2-5 \mu +6+\frac{4 \mu ^5+16 \mu ^4-7 \mu ^2-40 \mu -24}{4 (1+\mu )^{5/2}}\right)\\
& \hskip -1.cm \hbox{and}\ \mu(\tau,u)=\frac{(1-u)^2 }{u^2}\left(1+\frac{\tau}{3}\right)-1\ .
\end{split}
\label{eq:Laplaceperim}
\end{equation}
Introducing the inverse Laplace transform of $M(\mu)$, i.e.\ the quantity $\check{M}(X)$ such that 
$\displaystyle{\int_0^\infty e^{-\mu\, X}\check{M}(X)= M(\mu)}$, this quantity is easily computed and reads
\begin{equation*}
\hskip -1.cm \check{M}(X)= \frac{e^{-X}}{2\sqrt{\pi}} \left(-2 \sqrt{X} ((X-10) X-2)+e^X \sqrt{\pi } X (X (2 X-5)+6) \left(1-\text{erf}\left(\sqrt{X}\right)\right)\right)\ .
\end{equation*}
From the linear relation
\begin{equation}
\mu(\tau,u)=\frac{B(u)}{3}\, \tau+\big(B(u)-1\big)\ , \qquad  B(u)=\frac{(1-u)^2 }{u^2}\ ,
\label{eq:muBu}
\end{equation}
we immediately deduce, taking the inverse Laplace transform of \eqref{eq:Laplaceperim}, that
\begin{equation*}
\hskip -1.cm  D^{\rm out}(L,u)= \frac{(1-u)^6 }{u^3} \frac{3}{B(u)}\,  e^{(1-B(u)) X(u)} \check{M}(X(u))\ \hbox{with}\ X(u)=\frac{3 L}{B(u)}\ .
\end{equation*}
This is precisely the expression \eqref{eq:probainoutL} for $c=1/3$.
\vskip .3cm
To compute its counterpart $D^{\rm in}(L,u)$ in the in-regime, we use the explicit expression \eqref{eq:HkoneT} of $H(k,1,T)$ to derive the identity
\begin{equation*}
\begin{split}
&\hskip -1.cm H(k-k\, u,1,T(y))\Big|_{\rm sing.}=\frac{32 \, (\alpha_0+2 k (1-u)) \, \alpha_3}{\left(4 k^2 (1-u)^2+4 \alpha_0 \, k (1-u)+\alpha_0^2-1\right)^2} (1-12 h)^{3/2}\, \\
&\hskip -1.cm\hbox{whenever}\ Y\big(T(y)\big)= \alpha_0+\alpha_2 \eta^4+\alpha_3 \eta^6+\cdots \ ,\\ 
\end{split}
\end{equation*}
where, as before, $\eta=(1-12h)^{1/4}$.  It implies that
\begin{equation}
\begin{split}
& \hskip -1.cm \lim_{k\to \infty} E_{k,k\, u}^{\rm in}\left(e^{-\tau\,  L(k\, u)}\right)\\
&=\lim_{k\to \infty}
\frac{1}{{\frac{4}{7} k^3}}\Bigg(\frac{32 (\alpha_0(\tau,k\, u)+2 k (1-u)) \alpha_3(\tau,k\, u)}{\left(4 k^2 (1-u)^2+4 \alpha_0(\tau,k\, u)\, k  (1-u)+(\alpha_0(\tau,k\, u))^2-1\right)^2} \\
&\hskip 2cm - \frac{32 (\tilde{\alpha}_0(\tau,k\, u)+2 k (1-u)) \tilde{\alpha}_3(\tau,k\, u)}{\left(4 k^2 (1-u)^2+4 \tilde{\alpha}_0(\tau,k\, u)\, k  (1-u)+(\tilde{\alpha}_0(\tau,k\, u))^2-1\right)^2} \Bigg)\\
&\hskip -1.cm\hbox{whenever}\ Y\left(e^{-2 \frac{\tau}{(k\, u)^2}} T_{k\, u}(y)\right)= \alpha_0(\tau,k\, u)+\alpha_2(\tau,k\, u) \eta^4+\alpha_3(\tau,k\, u) \eta^6+\cdots \ ,\\
&\hskip -1.cm\hbox{and}\ Y\left(e^{-2 \frac{\tau}{(k\, u)^2}} T_{k\, u-1}(y)\right)= \tilde{\alpha}_0(\tau,k\, u)+\tilde{\alpha}_2(\tau,k\, u) \eta^4+\tilde{\alpha}_3(\tau,k\, u) \eta^6+\cdots \ .\\
\end{split}
\label{eq:Laplaceperimbis}
\end{equation}
Using the easily computed large $k$ expansions
\begin{equation*}
\begin{split}
&\hskip -1.cm \alpha_0(\tau,k\, u)= 2 \sqrt{\frac{3}{3+\tau}}\, (k\, u)+ 3 \left(\frac{3}{3+\tau}\right)^{3/2}+O\left(\frac{1}{k\, u}\right)\ ,\\
&\hskip -1.cm \tilde{\alpha}_0(\tau,k\, u)= 2 \sqrt{\frac{3}{3+\tau}}\, (k\, u)+  \left(\frac{3}{3+\tau}\right)^{3/2}+O\left(\frac{1}{k}\right)\ ,\\
&\hskip -1.cm \alpha_3(\tau,k\, u)=  \left(\frac{3}{3+\tau}\right)^{3/2} \, \frac{(k\, u)^7}{28}+ 3 \, (21+4\tau) \left(\frac{3}{3+\tau}\right)^{5/2} \, \frac{(k\, u)^6}{168}+O\left((k\, u)^5\right)\ ,\\
&\hskip -1.cm \tilde{\alpha}_3(\tau,k\, u)=  \left(\frac{3}{3+\tau}\right)^{3/2} \, \frac{(k\, u)^7}{28}+ (21+4\tau) \left(\frac{3}{3+\tau}\right)^{5/2} \, \frac{(k\, u)^6}{168}+O\left((k\, u)^5\right)\ ,\\
\end{split}
\end{equation*}
we deduce\footnote{Note that, for $\tau=0$, $ u^3 Q(\mu(0,u),B(u))=\frac{1}{4} (7-3 u)\, u^6=p^{\rm in}(u)$ as it should.}
\begin{equation*}
\begin{split}
& \hskip -1.cm \lim_{k\to \infty} E_{k,k\, u}^{\rm in}\left(e^{-\tau\,  L(k\, u)}\right)= u^3 Q(\mu(\tau,u),B(u))\\
& \hskip -1.cm \hbox{with}\ Q(\mu,B)=\frac{1}{\mu^4}
\left(-\!3 \mu ^2\!-\!3 B \mu\!-\!4 \mu\!-\!6 B\!+\!\frac{4 \mu ^3\!+\!3 B \mu ^2\!+\!20 \mu ^2\!+\!24 B \mu\!+\!16 \mu\!+\!24 B}{4 \sqrt{1+\mu}}\right)
\\
\end{split}
\end{equation*}
and $\mu(\tau,u)$ and $B(u)$ as in \eqref{eq:muBu}.
Introducing the inverse Laplace transform $\check{Q}(X,B)$ of $Q(\mu,B)$, easily computed to be
\begin{equation*}
\hskip -1.cm \check{Q}(X,B)=\frac{e^{-X}}{2\sqrt{\pi}} (B X+2) \left(2 \sqrt{X} (X+1)-e^X \sqrt{\pi } X (2 X+3) \left(1-\text{erf}\left(\sqrt{X}\right)\right)\right)\ ,
\end{equation*}
we immediately deduce, that
\begin{equation*}
\hskip -1.cm  D^{\rm in}(L,u)= u^3 \frac{3}{B(u)}\,  e^{(1-B(u)) X(u)} \check{Q}(X(u),B(u))\ \hbox{with}\ X(u)=\frac{3 L}{B(u)}\ .
\end{equation*}
This is precisely the expression \eqref{eq:probainoutL} for $c=1/3$.

\subsection{Joint law for the volume and perimeter at large ${\boldsymbol k}$ and ${\boldsymbol d}$ in the out-regime}
\label{sec:expljointlaw}
We now wish to control, in addition to the perimeter, the volume of the hull. Of course, this is non-trivial only if we condition the map configurations to be in the out-regime
where this volume is finite.  We are now interested in $G(k,d,g,h,\alpha)$ for arbitrary $\alpha$
and arbitrary $h=\rho/12$ (with $0\leq \rho\leq 1$). We consider again only the large $k$ and $d$ limit with fixed  $u=d/k$, a limit where ${\mathcal L}(d)$ growths like $d^2=(k\, u)^2$ 
while ${\mathcal V}(d)$ growths like $d^4=(k\, u)^4$. We therefore set
\begin{equation*}
\alpha=e^{-\frac{\tau}{(k\, u)^2}}\ , \qquad \rho=e^{-\frac{\sigma}{(k\, u)^4}}
\end{equation*}
with $\rho$ and $\sigma$ remaining finite.
From the relation \eqref{eq:gh} between $h$ and $y$,  taking the form of $\rho$ above amounts to setting
\begin{equation*}
y=e^{-\sqrt{6}\frac{\sigma^{1/4}}{k\, u}+O\left(\frac{1}{(k\, u)^3}\right)}\ .
\end{equation*}
We then have the expansions
\begin{equation*}
\begin{split}
&\hskip -.8cmY\left(\alpha^2 T_{k\,u}(y)\right)= 2\frac{1-u}{\sqrt{1+\mu}}\, k +9\sqrt{6}\, \frac{(1-u)^3}{u^3\, (1+\mu)^{3/2}} \frac{(1+W)(2+W)\sigma^{3/4}}{W^3}+O\left(\frac{1}{k}\right) \\
&\hskip -.8cmY\left(\alpha^2 T_{k\,u-1}(y)\right)= 2\frac{1-u}{\sqrt{1+\mu}}\, k +3\sqrt{6}\, \frac{(1-u)^3}{u^3\, (1+\mu)^{3/2}} \frac{(1+W)(2+W)\sigma^{3/4}}{W^3}+O\left(\frac{1}{k}\right) \\
&\hskip -.8cm \hbox{with}\ \mu\equiv \mu(\sigma,\tau,u)=\frac{(1-u)^2}{u^2}\, \left(
\frac{\tau}{3}+\frac{3\sqrt{\sigma }}{2}\left(\coth^2\left(\sqrt{\frac{3}{2}}\, \sigma^{1/4}\right)-\frac{2}{3}\right)\right)-1\\
&\hskip -.8cm\hbox{and}\ W\equiv W(\sigma)=e^{\sqrt{6} \, \sigma^{1/4}}-1\\
\end{split}
\end{equation*}
so that, eventually,
\begin{equation}
\begin{split}
& \hskip -1.cm \lim_{k\to \infty} E_{k,k\, u}^{\rm out}\left(e^{-\tau\,  L(k\, u)-\sigma\, V(k\, u)}\right)\\
&=\lim_{k\to \infty}\frac{1}{\frac{4}{7} k^3}\Bigg(
\mathfrak{h}_3\left(k\!-\!k\, u,
2\frac{1-u}{\sqrt{1+\mu}}\, k +9\sqrt{6}\, \frac{(1-u)^3}{u^3\, (1+\mu)^{3/2}} \frac{(1+W)(2+W)\sigma^{3/4}}{W^3} \right)\\
&\hskip 1.5cm \!-\!\mathfrak{h}_3\left(k\!-\!k\, u,
2\frac{1-u}{\sqrt{1+\mu}}\, k +3\sqrt{6}\, \frac{(1-u)^3}{u^3\, (1+\mu)^{3/2}} \frac{(1+W)(2+W)\sigma^{3/4}}{W^3}\right)\Bigg) \\
&=3 \sqrt{6}\ \frac{(u-1)^6}{u^3} \frac{(W+1) (W+2)}{W^3} \sigma ^{3/4}\, M(\mu)\\
&= \left(\frac{3}{2}\right)^{3/2}\ \frac{(u-1)^6}{u^3}\frac{ \cosh \left(\sqrt{\frac{3}{2}}\, \sigma^{1/4}\right)}{\sinh^3 \left(\sqrt{\frac{3}{2}}\, \sigma^{1/4}\right)}\, \sigma ^{3/4}\, M(\mu)
\end{split}
\label{eq:tausigmader}
\end{equation}
with $\mu=\mu(\sigma,\tau,u)$ and $W=W(\sigma)$ as above, and where the function $M(\mu)$ has the same expression as in \eqref{eq:Laplaceperim}.
Normalizing by $p^{\rm out}(u)$, this yields precisely the announced expression \eqref{eq:expsigmatau} with $c=1/3$ and $f=36$.
 
From the linear relation between $\mu$ and $\tau$
\begin{equation*}
\begin{split}
&\hskip -1.cm \mu(\sigma,\tau,u)= B(u)\, \frac{\tau}{3}+\big(A(\sigma) B(u)-1\big)\ , \\
&\hskip -1.cm \hbox{with}\  B(u)=\frac{(1-u)^2 }{u^2}\ \  \hbox{and}\  \ A(\sigma)=\frac{3\sqrt{\sigma }}{2}\left(\coth^2\left(\sqrt{\frac{3}{2}}\, \sigma^{1/4}\right)-\frac{2}{3}\right)\ , \\
\end{split}
\end{equation*}
we immediately deduce, taking the inverse Laplace transform of \eqref{eq:tausigmader}, that
\begin{equation*}
\begin{split}
&\hskip -1.2cm \lim_{k\to \infty} E_{k}\left(e^{-\sigma\, V(k\, u)}\Big|{\mathcal V}(k\, u)\ \hbox{finite and}\ L(k\, u)
 =L\right) \\
 &= \frac{1}{D^{\rm out}(L,u)}\left(\frac{3}{2}\right)^{3/2}\ \frac{(u-1)^6}{u^3}\frac{ \cosh \left(\sqrt{\frac{3}{2}}\, \sigma^{1/4}\right)}{\sinh^3 \left(\sqrt{\frac{3}{2}}\, \sigma^{1/4}\right)}\, \sigma ^{3/4}\\
 & \hskip 2.cm \times  \frac{3}{B(u)}\,  e^{(1-A(\sigma)\, B(u)) X(u)} \check{M}(X(u))\quad  \hbox{with}\ X(u)=\frac{3 L}{B(u)}\\
& =\left(\frac{3}{2}\right)^{3/2} \frac{ \cosh \left(\sqrt{\frac{3}{2}}\, \sigma^{1/4}\right)}{\sinh^3 \left(\sqrt{\frac{3}{2}}\, \sigma^{1/4}\right)}\, \sigma ^{3/4}\, e^{\left(1-A(\sigma)\right)\, B(u) X(u)}\\
& =\left(\frac{3}{2}\right)^{3/2} \frac{ \cosh \left(\sqrt{\frac{3}{2}}\, \sigma^{1/4}\right)}{\sinh^3 \left(\sqrt{\frac{3}{2}}\, \sigma^{1/4}\right)}\, \sigma ^{3/4}\, e^{-3\, L\, \left(A(\sigma)-1\right)}\ .\\
   \end{split}
\end{equation*}
This is precisely the announced result \eqref{eq:VcondL} for $c=1/3$ and $f=36$. Remarkably, all the $u$ dependences dropped out upon normalizing by $D^{\rm out}(L,u)$ and
the above conditional probability is thus \emph{independent of $u$}.
 
\section{Other families of maps}
\label{sec:Other}
The strategy presented in this paper may be applied to other families of maps provided that a coding by slices exists and that the decomposition of the corresponding
slices along dividing lines at a fixed distance from their apex is fully understood. This is the case for planar triangulations, as explained in \cite{G15a}, and
for planar Eulerian triangulations, as explained in \cite{G16b}.
We have reproduced and adapted the computations above to deal with these two other families of maps. We do not display here our calculations 
since they are quite tedious and give no really new information. Indeed, we find that, in the limit
of large $k$ and $d$ with fixed $u=d/k$, all the laws that we obtained for quadrangulations have exactly the same expressions for triangulations and Eulerian triangulations, up
to a global normalization for the rescaled length $L(d)=\mathcal{L}(d)/d^2$ and a global normalization for the rescaled volume $V(d)=\mathcal{V}(d)/d^4$.
If we adopt the definitions of \cite{G15a} and \cite{G16b} for the hull at distance $d$ in triangulation and Eulerian triangulations respectively, 
these normalizations amount to change in our various laws of Section~\ref{sec:summary} the values $c=1/3$ and $f=36$ found in Section~\ref{sec:explicit} for quadrangulations 
to the values displayed in \eqref{eq:cval} and \eqref{eq:fval}.

The origin of the scaling factor $f$ is easily found in the relation between the weight $h$ per face in the hull and the variable $y$
which is ``conjugated" to the distance $d$ in the slice generating function $T_d(y)$ (by this, we mean that $d$ appears in $T_d(y)$ via 
the combination $y^d$ only). We have for the three families (i), (ii) and (iii) of maps (see for instance \cite{G15a,G15b,G16b})
\begin{equation*}
\hskip -1.cm \hbox{(i):}\  h(y)=\frac{y \left(1+y+y^2\right)}{\left(1+4y+y^2\right)^2}\ , \ \ 
\hbox{(ii):}\  h(y)= \frac{\sqrt{y (1+y)}}{\left(1+10 y+ y^2\right)^{3/4}}\ , \ \ 
\hbox{(iii):}\  h(y)= \frac{y \left(1+y^2\right)}{(1+y)^4}\ .
\end{equation*}
The desired singularities are obtained for $h(y)=h(1)\rho$ where $\rho=e^{-\sigma/d^4}$ if we wish to capture 
the large $d$ behavior of the rescaled volume. Setting $h(y)=h(1)e^{-\sigma/d^4}$ amounts, at large $d$ to setting
$y=e^{-(f\, \sigma)^{1/4}/d}$ with, from the above relations $f=36$ in case (i), $f=192$ in case (ii) and $f=16$ in case (iii).
All the universal laws involving $y^d\simeq e^{-(f\sigma)^{1/4}}$, the quantity $\sigma$ always appears via the combination
$(f \, \sigma)^{1/4}$ in our various laws. This explains the origin of $f$.
\vskip .3cm
To understand the origin of the normalization factor $c$, the simplest quantity to compute is probably
the limiting expectation value
\begin{equation*}
\lim_{k\to\infty}E_{k,d}^{\rm out}\left(\alpha^{\mathcal{L}(d)}\right)\ .
\end{equation*}
In the case of quadrangulations, we have from \eqref{eq:Hsing} the large $k$ expansion
\begin{equation*}
\mathfrak{h}_3(k,T)=\frac{1}{7}\, k^4 + \frac{2\, Y}{7}\, k^3 + O(k^2)
\end{equation*}
so that (after normalization by $\mathfrak{f}_3(k)\sim (4/7)k^3$)
\begin{equation*}
\begin{split}
&\lim_{k\to\infty}E_{k,d}^{\rm out}\left(\alpha^{\mathcal{L}(d)}\right)=\frac{1}{2}\left(Y\left(\alpha^2\, T_d(1)\right)-Y\left(\alpha^2\, T_{d-1}(1)\right)\right)\\
&
= \frac{1}{2} \left(\sqrt{\frac{6 \alpha ^2+(d+1) (d+2) \left(25-\alpha ^2\right)}{6 \alpha ^2+(d+1) (d+2) \left(1-\alpha
   ^2\right)}}-\sqrt{\frac{6 \alpha ^2+d (d+1) \left(25-\alpha ^2\right)}{6 \alpha ^2+d (d+1) \left(1-\alpha ^2\right)}}\right)\ .\\
  \end{split}
\end{equation*}
Using (for $C>0$)
\begin{equation*}
\hskip -1.2cm \sqrt{\frac{C^2-\beta}{1-\beta}}=C+2\sum_{p\geq 1}\beta^p\, A_p(C)\ \ \hbox{with}\ A_p(C)=\frac{1}{C^{2p-1}}\, \sum_{q=0}^{p-1}{p-1\choose q}{2q+1\choose q}\left(\frac{C^2-1}{4}\right)^{q+1}
\end{equation*}
here with $C^2=25$, we deduce that
\begin{equation*}
\hbox{(i):}\ \lim_{k\to\infty}P_{k,d}^{\rm out}\left({\mathcal{L}(d)=2p}\right)=A_p(5)\, \left(\left(\frac{(d-1)(d+4)}{d+1)(d+2)}\right)^p-\left(\frac{(d-2)(d+3)}{d(d+1)}\right)^p\right)\ ,
\end{equation*}
where the subscript ``out" is irrelevant since for finite $d$ and infinite $k$, map configurations are necessarily in the out-regime.
A similar calculation for the families (ii) and (iii) yields
\begin{equation*}
\begin{split}
& \hbox{(ii):}\ \lim_{k\to\infty}P_{k,d}^{\rm out}\left({\mathcal{L}(d)=p}\right)=A_p(3)\, \left(\left(\frac{d(d+3)}{d+1)(d+2)}\right)^p-\left(\frac{(d-1)(d+2)}{d(d+1)}\right)^p\right)\ ,\\
& \hbox{(iii):}\ \lim_{k\to\infty}P_{k,d}^{\rm out}\left({\mathcal{L}(d)=2p}\right)=2 A_p(3)\, \left(\left(\frac{(d-1)(d+5)}{d+1)(d+3)}\right)^p-\left(\frac{(d-2)(d+4)}{d(d+2)}\right)^p\right)\\
\end{split}
\end{equation*}
(note that $\mathcal{L}(d)$ is necessarily even in case (iii) but has arbitrary parity in case (ii)).
From the large $p$ behavior $A_p(C)\sim \sqrt{C^2-1}/(2\sqrt{\pi\, p})$, we immediately deduce, taking $d$ and $p$ large with $p/d^2=L/2$ (case (i) and (iii)) or
 $p/d^2=L$ (case (ii)) the following probability densities for the three families of maps:
\begin{equation*}
\begin{split}
& \hbox{(i):}\  \lim_{d\to \infty} \left(\lim_{k\to\infty} \frac{1}{dL} P^{\rm out}_{k,d}\left(L\leq L(d) < L+dL\right)\right)= 6\, \sqrt{3}\, \sqrt{\frac{L}{\pi}}\, e^{-3\, L}\ , \\
& \hbox{(ii):}\  \lim_{d\to \infty} \left(\lim_{k\to\infty} \frac{1}{dL} P^{\rm out}_{k,d}\left(L\leq L(d) < L+dL\right)\right)= 4\, \sqrt{2}\, \sqrt{\frac{L}{\pi}}\, e^{-2\, L}\ , \\
& \hbox{(iiI):}\  \lim_{d\to \infty} \left(\lim_{k\to\infty} \frac{1}{dL} P^{\rm out}_{k,d}\left(L\leq L(d) < L+dL\right)\right)= 16\, \sqrt{\frac{L}{\pi}}\, e^{-4\, L}\ . \\
\end{split}
\end{equation*}
In agreement with the equivalence principle \eqref{eq:utozero}, these law reproduce the general form \eqref{eq:uzeroone}
for the limit $u\to 0$ of $D^{\rm out}(L,u)/p^{\rm out}(u)$ (recall that $p^{\rm out}(u)
\to 1$ for $u\to 0$) via the identification $c=1/3$ in case (i), $c=1/2$ in case (ii) and $c=1/4$ in case (iii).

We end this section by giving for bookkeeping the (non-universal) expression for the probability $P_k({\mathcal V}(d)\ \hbox{infinite})$ at finite $k$ and $d$ for the families (ii) and (iii).
We find
\begin{equation*}
\begin{split}
&\hskip -1.cm \hbox{(ii):}\  P_k({\mathcal V}(d)\ \hbox{infinite})=\frac{k^2 (k+1)^2}{2 (2 k+1) \left(5 k^6+15 k^5+14 k^4+3 k^3-k^2-1\right)}\\ 
& \hskip 1.cm \times \Bigg(\frac{d (d+1) (d+2) (d+3) \left(10 d^4+60 d^3+146 d^2+168 d+71\right)}{(2 d+3) (k+1)^3}\\
&\hskip 2.cm  - \frac{(d-1) d (d+1) (d+2) \left(10 d^4+20 d^3+26 d^2+16 d-1\right)}{(2 d+1) k^3}\Bigg) \ ,\\
&\hskip -1.cm \hbox{(iii):}\  P_k({\mathcal V}(d)\ \hbox{infinite})= \frac{k (k+1) (k+2) (k+3)}{2 (2 k+3) \left(10 k^6+90 k^5+283 k^4+348 k^3+103 k^2-42 k-36\right)}\\ 
& \hskip 1.cm\times \Bigg( \frac{(d-1) (d+1) (d+3) (d+5) \left(10 d^4+80 d^3+256 d^2+384 d+189\right) (k+2)}{(d+2) (k+1)^2 (k+3)^2} \\
& \hskip 2.cm  - \frac{(d-2) d (d+2) (d+4) \left(10 d^4+40 d^3+76 d^2+72 d-9\right) (k+1)}{ (d+1) k^2 (k+2)^2}\Bigg)\ .\\
\end{split}
\end{equation*}
When $k,d\to \infty$ and $d/k=u$, both expressions tend to $p^{\rm in}(u)=(7-3u)u^6/4$.
 
\section{Conclusion}
\label{sec:conclusion}
\begin{figure}
\begin{center}
\includegraphics[width=8cm]{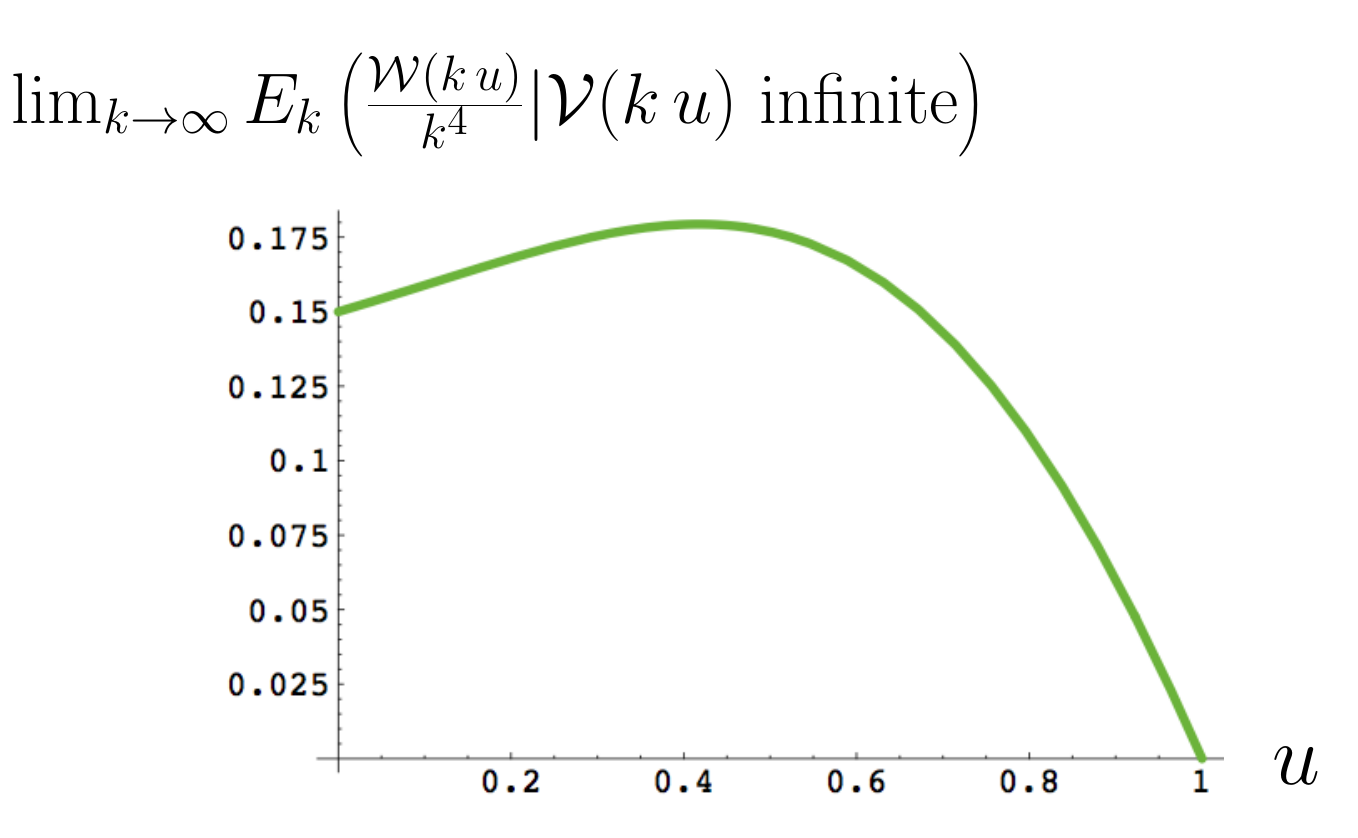}
\end{center}
\caption{A plot of the expectation value of the (properly normalized) volume $\mathcal{W}(d)\equiv N-\mathcal{V}(d)$ in the in-regime as a function of $u=k/d$ for large $k$ and $d$
(here with $f=36$).}
\label{fig:complementvol}
\end{figure} 

In this paper, we explored the statistics of hull perimeters for three families of infinitely large planar maps: quadrangulations, triangulations and Eulerian triangulations, 
with a particular emphasis on the influence on this statistics of the constraint that the map configurations either yield a finite hull volume 
or not. In the case where the hull volume is finite, we also discussed the statistics of this volume itself, as well as its coupling to the hull perimeter statistics.
Our study, based on an accurate coding of $k$-pointed-rooted planar maps by $k$-slices, makes a crucial use of a particular  
 \emph{recursive decomposition of these slices} obtained by cutting them along lines which precisely follow hull boundaries for increasing distances $d$ (see figure \ref{fig:constrhull}
 for an illustration) for $d<k$. This decomposition, initiated in \cite{G15a} for triangulations, and then extended in \cite{G15b,G16b} for the two other families of maps, may be used 
 to address many other questions of the type discussed here, either for the same geometry, i.e.\ within pointed-rooted maps, or for other more involved geometries.
 
 Among other quantities which may be computed within the above geometry of $k$-pointed-rooted maps are the statistics of the volume $\mathcal{W}(d)\equiv N-\mathcal{V}(d)$
of the complementary of the hull at distance $d$, i.e.\ the component containing the marked vertex $x_1$. In the in-regime, $\mathcal{W}(d)$
is finite and we may compute its limiting universal expectation value for large $k$ and $d$. We find
\begin{equation*}
\hskip -1.cm \lim_{k\to \infty}E_{k}\left(\frac{\mathcal{W}(k\, 
u)}{k^4}\Big|\mathcal{V}(k\, u)\ \hbox{infinite}\right)=f\, \frac{(1-u) \left(14+16\, u+16\, u^2+16\, u^3-39\, u^4+12\, u^5\right)}{480 (7-3 u)}\ .
\end{equation*}
This quantity is plotted in figure \ref{fig:complementvol} for $f=36$ (quadrangulations).

Concerning other tractable geometries, we recall that pointed \emph{maps with a boundary} (i.e.\ maps with a distinguished external face of arbitrary degree) may be decomposed into sequences of slices and our recursive decomposition of slices gives a direct access to the statistics of a generalized hull at distance $d$ whose boundary would separate the pointed vertex from the external face (assuming that all vertices of the boundary are at a distance strictly larger that $d$ from the pointed vertex). 
  
To conclude, many other families of maps (for instance maps with prescribed face degrees) may be coded by slices and, even if a recursion relation
of the type of Ref.~\cite{G15a,G15b,G16b} is not known in general for these slices\footnote{Other recursions are known however.}, the actual
form of the associated slice generating functions is known in many cases \cite{BG12}. This might be enough to address the hull statistics for these maps since,
as the reader noticed, the actual expression for the operator $\mathcal{K}$ describing the action of one step of the recursion is
not really needed. What is needed is an equation of the form \eqref{eq:propK} which displays the result of this operator on properly parametrized generating functions.
This equation itself is moreover directly read off the explicit expression of the slice generating functions themselves (here for quadrangulations).
Slices associated with maps with arbitrary face degrees have generating functions whose expressions are of the same general form 
(although more involved in general) as that for quadrangulations (see \cite{BG12}). The actual knowledge of these expressions might thus be sufficient 
to infer the hull statistics for the corresponding maps.

\appendix
\section{Expectation value of the perimeter at finite $k$ and $d$ in the out- and in-regimes} 
Computing the expectation value of the perimeter simply involves computing the quantity $\partial_\alpha G(k,d,g,h,\alpha)\Big|_{\alpha=1}$, which 
itself, from \eqref{eq:Gexpr}, simply requires an expression
for the quantity
\begin{equation*}
2 T\, \partial_T H(k,x,T)\ .
\end{equation*}
In the out-regime, we need to estimate the singularity of this latter quantity when $g\to 1/12$ ($x\to 1$). We find
\begin{equation*}
\begin{split}
&
\hskip -1.cm 2 T\, \partial_T H(k,x,T)\Big|_{\rm sing.}=2 T\, \partial_T \mathfrak{h}_3\big(k,Y(T)\big) (1-12 g)^{3/2}\\
&\hskip 2.cm = \mathfrak{dh}_3\big(k,Y(T)\big) (1-12 g)^{3/2}\\
& \hskip -1.cm \hbox{with}\  \mathfrak{dh}_3(k,Y)=\frac{(25-Y^2)(1-Y^2)}{24Y}\partial_Y\mathfrak{h}_3(k,Y)
\\ &= \frac{k (25-Y^2)(1-Y^2)}{20160 Y^3 (2 k+Y-1)^3 (2 k+Y+1)^3} \big(315 Y^{10}+3780 k Y^9+19740 k^2 Y^8-1995 Y^8\\ & +60480 k^3 Y^7-20160 k Y^7+120960 k^4 Y^6-82320 k^2 Y^6+16590 Y^6+161280 k^5
   Y^5\\ & -174720 k^3 Y^5+71400 k Y^5+138240 k^6 Y^4-209664 k^4 Y^4+101640 k^2 Y^4+3594 Y^4\\ & +69120 k^7 Y^3-139776 k^5 Y^3+60480 k^3 Y^3-26784 k Y^3+15360 k^8 Y^2-39936 k^6
   Y^2\\ & +24192 k^4 Y^2-54224 k^2 Y^2-36217 Y^2-65100 k Y-21700 k^2+5425\big)\ .\\
\end{split}
\end{equation*}
Setting $h=1/12$ ($y=1$) and $\alpha=1$ so that the values of interest are $Y(T_d(1))=2d+3$ and $Y(T_{d-1}(1))=2d+1$,  we immediately deduce, upon normalization, that
\begin{equation*}
E_{k}\left({\mathcal L}(d)\Big|{\mathcal V}(d)\ \hbox{finite}\right)=\frac{\mathfrak{dh}_3(k-d,2d+3)-\mathfrak{dh}_3(k-d,2d+1)}{\mathfrak{h}_3(k-d,2d+3)-\mathfrak{h}_3(k-d,2d+1)}\ .
\end{equation*}
This immediately yields an explicit expression (which we do not reproduce here) for the expectation value of the perimeter at finite $k$ and $d$ in the out-regime. It is then easily verified that at large $k$ and $d$, ${\mathcal L}(d)$ scales as $d^2$ and that, 
for $k,d\to \infty$ and $d/k=u$ fixed, the expression for the expectation value of $L(d)={\mathcal L}(d)/d^2$ simplifies into the formula given in the first line of \eqref{eq:expectperim},
with here $c=1/3$.

In the in-regime, we now set $g=1/12$ ($x=1$). Using, from \eqref{eq:HkoneT},
\begin{equation*}
2T\, \partial_T H(k,1,T)=\frac{\big(25-Y(T)^2\big)\big(1-Y(T)^2\big)}{24Y(T)} \frac{32 (2 k+Y(T))}{(2 k+Y(T)-1)^2 (2 k+Y(T)+1)^2}
\end{equation*}
and plugging the expansion \eqref{eq:expTd} for $Y(T_d(y))$ when $y\to 1$ ($\eta\to 0$), we deduce
\begin{equation*}
\begin{split}
&
\hskip -1.cm 2T\, \partial_T H\big(k,1,T_d(y)\big)\Big|_{\rm sing.}=\widetilde{\delta\mathfrak{h}}_3(k,d) (1-12 h)^{3/2}\\
&\hskip -1.cm \hbox{with}\ \widetilde{\delta\mathfrak{h}}_3(k,d)=\frac{2 (d-1) (d+1) (d+2) (d+4)}{315 (2 d+3)^3 (d+k+1)^3
   (d+k+2)^3} \left(15 d^4+90 d^3+237 d^2+306 d+140\right)\\ & \times \big(6 k d^6+12 k^2 d^5+54 k d^5+24 d^5+6 k^3 d^4+90 k^2 d^4+240 k d^4+180 d^4+36 k^3
   d^3\\ & +270 k^2 d^3+630 k d^3+534 d^3+68 k^3 d^2+405 k^2 d^2+937 k d^2+783 d^2+42 k^3 d\\ &+285 k^2 d+705 k d+567 d-2 k^3+63 k^2+203 k+162\big)\ .
   \\
   \end{split}
\end{equation*}
We immediately deduce, upon normalization, that
\begin{equation*}
E_{k}\left({\mathcal L}(d)\Big|{\mathcal V}(d)\ \hbox{infinite}\right)=\frac{\widetilde{\mathfrak{dh}}_3(k-d,d)-\widetilde{\mathfrak{dh}}_3(k-d,d-1)}{\tilde{\mathfrak{h}}_3(k-d,d)-\tilde{\mathfrak{h}}_3(k-d,d-1)}
\end{equation*}
which again yields an explicit expression (not reproduced here) for the expectation value of the perimeter at finite $k$ and $d$ in the in-regime.
It is again easily verified that, 
for $k,d\to \infty$ and $d/k=u$ fixed, the expression for the expectation value of $L(d)={\mathcal L}(d)/d^2$ simplifies into the formula given in the second line of \eqref{eq:expectperim},
with here $c=1/3$.

\begin{figure}
\begin{center}
\includegraphics[width=8cm]{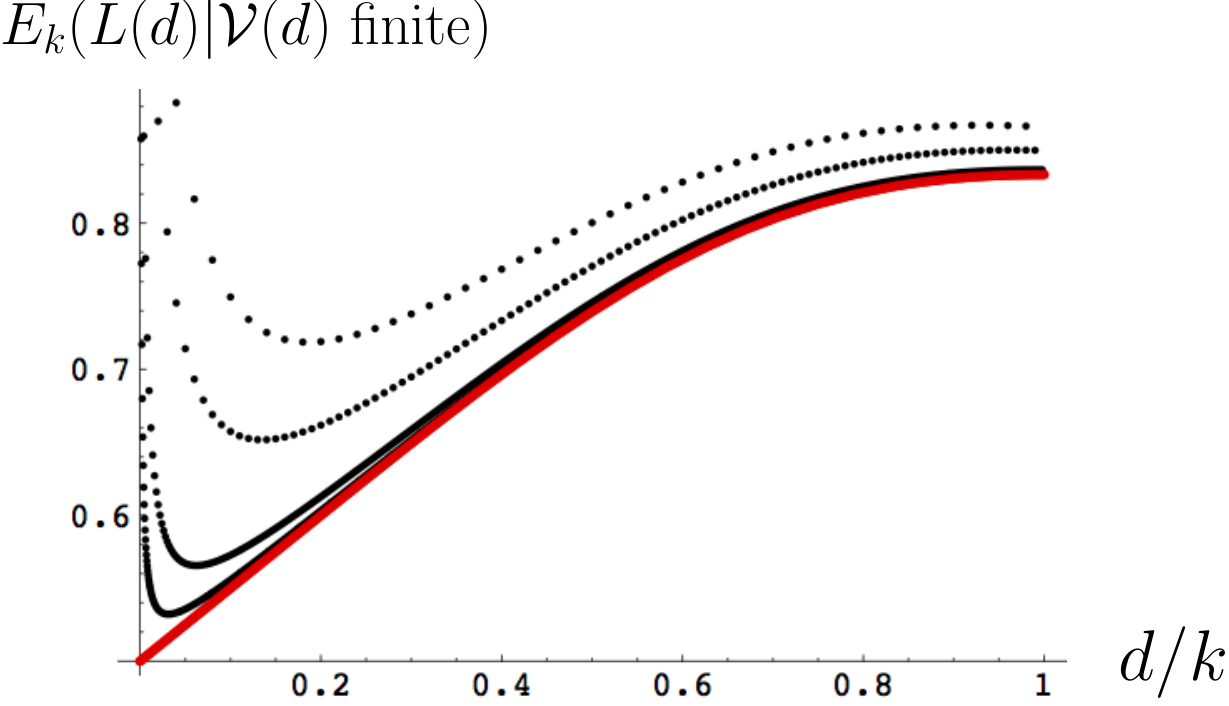}
\end{center}
\caption{Plots of the expectation value of the rescaled hull perimeter $L(d)$ in the out-regime as a function of $d/k$ for $k=50$, $100$, $500$, and $2000$. In red:
the corresponding limiting law for large $k$ and $d$, as given by the first line of \eqref{eq:expectperim}.}
\label{fig:perimeteroutfinite}
\end{figure}

\begin{figure}
\begin{center}
\includegraphics[width=8cm]{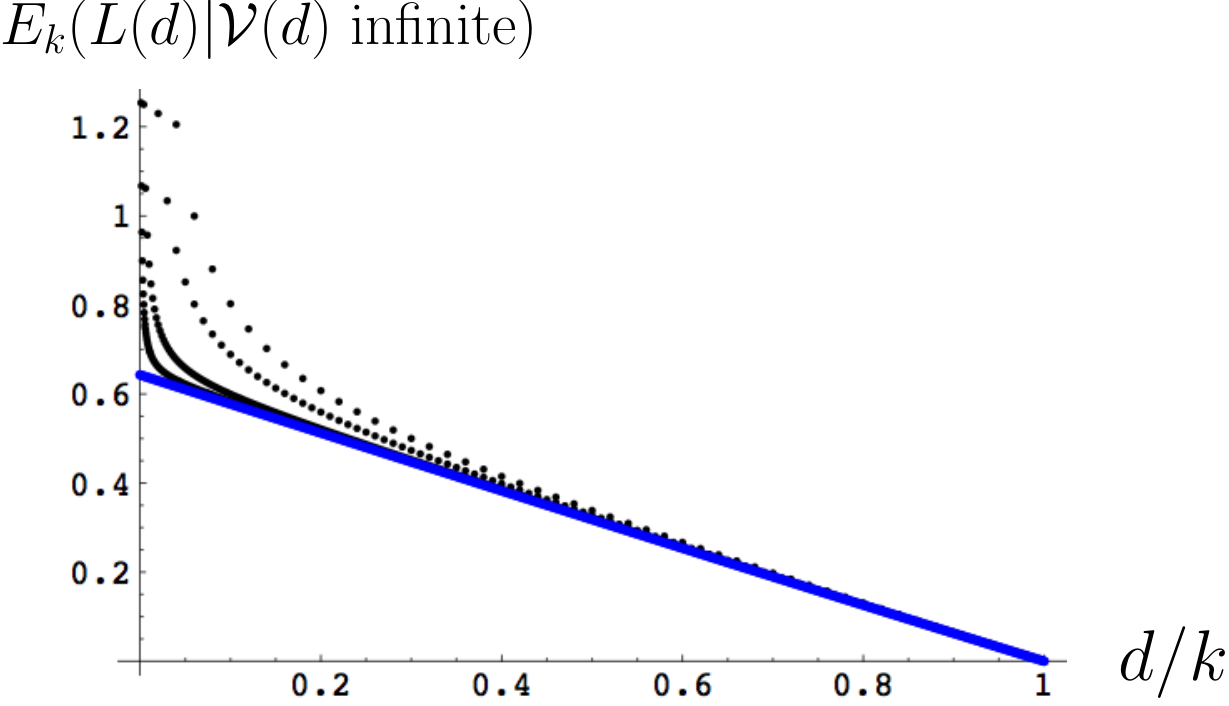}
\end{center}
\caption{Plots of the expectation value of the rescaled hull perimeter $L(d)$ in the in-regime as a function of $d/k$ for $k=50$, $100$, $500$, and $2000$. In blue:
the corresponding limiting law for large $k$ and $d$, as given by the second line of \eqref{eq:expectperim}.}
\label{fig:perimeterinfinite}
\end{figure}

Figure~\ref{fig:perimeteroutfinite} (respectively figure~\ref{fig:perimeterinfinite}) shows a comparison 
between the limiting expression given in the first (respectively the second) line of \eqref{eq:expectperim}
with $c=1/3$ vs $u$ and the finite $k$ and $d$ expression for
$E_{k}\left(L(d)\Big|{\mathcal V}(d)\ \hbox{finite}\right)$ (respectively $E_{k}\left(L(d)\Big|{\mathcal V}(d)\ \hbox{infinite}\right)$) vs $d/k$ for $k=50$, $100$, $500$, and $2000$ 
and $2\leq d\leq k-1$.

\section*{Acknowledgements} 
The author acknowledges the support of the grant ANR-14-CE25-0014 (ANR GRAAL).
 
\bibliographystyle{plain}
\bibliography{hullrefined}

\end{document}